\documentclass[10pt]{amsart}
 \usepackage{amssymb,latexsym,amsmath,amsthm}
  \usepackage{fullpage}

 \renewcommand{\div}{\mathop{\mathrm{div}}\nolimits}

 \usepackage[numbers,sort&compress]{natbib}

\newtheorem{thm}{Theorem}[section]

\newtheorem{dfn}{Definition}[section]

\newtheorem{lemma}{Lemma}[section]
\newtheorem{notation}{Notation}[section]

\newtheorem{rem}{Remark}[section]

\numberwithin{equation}{section}

\begin{document}
\date{}

\title{Symmetry results for fractional elliptic systems and related problems}

\author{Mostafa Fazly}
\author{Yannick Sire}

\address{Department of Mathematical and Statistical Sciences, University of Alberta, Edmonton, Alberta, Canada T6G 2G1}

\email{fazly@ualberta.ca}

\address{Universite Aix-Marseille and LATP 9, rue F. Joliot Curie, 13453 Mar- seille Cedex 13, France.}

\email{sire@cmi.univ-mrs.fr}

\thanks{The first author  is pleased to acknowledge the support of a University of Alberta start-up grant RES0019810. The second author is supported by the ANR project "HAB" and the ERC Starting grant "EPSILON"}

\maketitle

\tableofcontents

\begin{abstract}
We study elliptic gradient systems with fractional laplacian operators on the whole space 
$$ (- \Delta)^\mathbf s \mathbf u =\nabla H (\mathbf u) \ \ \text{in}\ \ \mathbf{R}^n,$$  
where  $\mathbf u:\mathbf{R}^n\to \mathbf{R}^m$, $H\in C^{2,\gamma}(\mathbf{R}^m)$ for $\gamma > \max(0,1-2\min \left \{s_i \right \})$, $\mathbf s=(s_1,\cdots,s_m)$ for $0<s_i<1$ and $\nabla H (\mathbf u)=(H_{u_i}(u_1, u_2,\cdots,u_m))_{i}$. We prove De Giorgi type results for this system for certain values of $\mathbf s$ and in lower dimensions, i.e. $n=2,3$.  Just like the local case, the concepts of orientable systems and $H-$monotone solutions, established in \cite{FG}, play the key role in proving symmetry results.   In addition, we provide optimal energy estimates, a monotonicity formula, a Hamiltonian identity and various Liouville theorems.  
\end{abstract}

\section{Introduction and main results}

This paper is devoted to several symmetry results, qualitative properties and Liouville-type theorems for solutions of non local elliptic systems.  Non local equations have led to several active research areas, in both applied and purely theoretical aspects for the past recent years. The prototypical operator involved is the so-called fractional laplacian $(-\Delta)^s$ for $s \in (0,1)$. It is a Fourier multiplier of symbol $|\xi|^{2s}$, see \cite{landkof}. Despite its interest in harmonic analysis (see the book \cite{landkof} for the first systematic study of potential analysis of this operator) it is also of great importance in the probability theory. Indeed, the fractional laplacian is the basic example of infinitesimal generator of Levy processes, see the book of Bertoin \cite{B} for an extensive study of such stochastic processes. Levy processes are processes  whose generators are given by the following formula (up to a normalizing constant) 
$$\mathcal I u(\mathbf x)=\int_{\mathbf R^n} (u(\mathbf x+\mathbf y)+u(\mathbf x-\mathbf y)-2u(\mathbf y))\mu(d\mathbf y)$$
for sufficiently smooth (say) functions $u$ and where $\mu(d\mathbf y)$ is a Levy measure, i.e. a positive function on $\mathbf R^n$ such that $\mu( \left \{ 0 \right \})=0$ and 
$$\int_{\mathbf R^n} \min (|\mathbf x|^2,1)\mu(d\mathbf x) < \infty.$$
In the case of the fractional laplacian, one has 
$$(-\Delta)^s u(\mathbf x)=\text{P.V.}\int_{\mathbf R^n} \frac{u(\mathbf x)-u(\mathbf y)}{|\mathbf x-\mathbf y|^{n+2s}}\,d\mathbf y,$$
where $\text{P.V.}$ stands for the principal value in the Cauchy sense.  

We consider the following system
\begin{eqnarray}
\label{main}
(- \Delta)^\mathbf s \mathbf u =\nabla H (\mathbf u) \ \ \text{in}\ \ \mathbf{R}^n,
  \end{eqnarray}
where $\mathbf u:\mathbf{R}^n\to \mathbf{R}^m$, $H\in C^{2,\gamma}(\mathbf{R}^m)$ for $\gamma > \max(0,1-2\min \left \{s_i \right \})$, $\mathbf s=(s_1,\cdots,s_m)$ where $0<s_i<1$ and $\nabla H (\mathbf u)=(H_{u_i}(u_1, u_2,...u_m))_{i}$. The notation $H_{u_i}$ stands for the partial derivative $\frac{\partial H}{\partial u_i}$. Therefore each component satisfies the equation 
$$ (- \Delta)^{s_i} u_i = \partial_{u_i}H(\mathbf u) \ \ \text{in}\ \ \mathbf{R}^n.
$$
For the local problems that is when $\mathbf s=1$, the above system has been studied for various purposes because of the interesting structure of the system. We refer interested readers to \cite{FG, ali} for symmetry results and to \cite{CF} for the regularity of extremal solutions of eigenvalue problems.

It is now a well known fact, and extensively used for these equations, that the fractional laplacian can be realized as the boundary operator (more precisely the Dirichlet-to-Neumann operator) of a suitable extension in the half-space (see \cite{cafS}). In view of this result, we will be considering the following extended system of equations 
 \begin{eqnarray}\label{emain}
 \left\{ \begin{array}{lcl}
\hfill \div(y^{a_i} \nabla v_i)&=& 0   \ \ \text{in}\ \ \mathbf{R}_+^{n+1}=\left \{x \in \mathbf R^n, y>0 \right \},\\   
\hfill -\lim_{y\to0}y^{a_i} \partial_{y} v_i&=& d_{s_i} \partial_{v_i}H(\mathbf v)   \ \ \text{in}\ \ \partial\mathbf{R}_+^{n+1},
\end{array}\right.
  \end{eqnarray}
where $a_i=1-2s_i$ and $d_{s_i} = \frac{\Gamma(1-s_i)}{2^{2s_i-1}\Gamma(s_i)}$.  Here $v_i$ is the extension of the function $u_i$.

The main results of the present paper deal with the solutions of \eqref{emain}. The direct consequences of these results, in the light of \cite{cafS}, leads us to similar results for the original system \eqref{main}. 

Our results are inspired by a famous conjecture of De Giorgi announced in \cite{DeG}. This conjecture concerns the flatness of level-sets of bounded monotone solutions of the scalar Allen-Cahn equation. The De Giorgi's conjecture is known to be true in $n=2$ by Ghoussoub-Gui \cite{GG1}, in $n=3$ by Ambrosio-Cabre \cite{AC}, in $4 \leq n \leq 8$ by Savin \cite{savin} (with an additional natural hypothesis). A counterexample is provided  in $n\ge 9$ by Del Pino-Kowalczyk-Wei \cite{PKW}. In addition, Fazly-Ghoussoub in  \cite{FG} established De Giorgi type results for elliptic systems of the form $\Delta \mathbf u=\nabla H(\mathbf u)$ where $\mathbf u:\mathbf R^n\to\mathbf R^m$ in dimensions $n=2,3$.    Corresponding symmetry results for non local equations are provided by Cabr\'{e}-Sire in \cite{CS2} and by Sire-Valdinoci in \cite{SV} when $n=2$ and  by Cabr\'{e}-Cinti in \cite{cabre} for $n=3$. Moreover, Dipierro-Pinamonti in \cite{serena} provided symmetry results for the system (\ref{main}) when $n=m=2$. 


 Before stating our main results, we would like to define the following concepts. 
\begin{dfn} \rm
We say that a solution ${\bf u}=(u_i)_{i=1}^m$ of (\ref{main}) is $H$-monotone if the following hold,
\begin{enumerate}
 \item For every $i\in \{1,\cdots, m\}$, $v_i$ is strictly monotone in the $x_n$-variable (i.e., $\partial_n v_i\neq 0$).

\item  For $i<j$, 
we have 
  \begin{equation}
\hbox{$\partial_{u_iu_j}H({\mathbf u}) \partial_n u_i(\mathbf x) \partial_n u_j (\mathbf x)\ge 0$  for all $x\in\mathbf{R}^n$.}
\end{equation}

\end{enumerate}
\end{dfn}

\begin{dfn} \label{weak} \rm 
 We shall say that the system  (\ref{main}) is orientable, if there exist nonzero functions $\theta_k\in C^1(\mathbf{R}^{n+1}_+)$, $k=1,\cdots,m$, which do not change sign, such that for all $i,j$ with $1\leq  i<j\leq m$, we have 
 \begin{equation}\label{oriantableu}
 \hbox{$ \partial_{u_iu_j}H({\mathbf u}) \theta_i(\mathbf x)\theta_j(\mathbf x)\ge 0$ \, for all $\mathbf x\in\mathbf{R}^n$.}
  \end{equation} 
  Similarly, if the condition (\ref{oriantableu}) holds for the extension $\mathbf v$, then we say (\ref{emain}) is orientable. 
\end{dfn}
This is a combinatorial assumption on the sign of the nonlinearity $H$.
 
\begin{dfn} \rm 
\begin{enumerate}
\item A solution $\mathbf v=(v_i)_i$ of the system (\ref{emain}) in $\mathbf{R}_+^{n+1}$ is said to be
{\it pointwise-stable}, if there exist $(\phi_i)_{i=1}^m$ in $C^\infty(\mathbf{R}_+^{n+1})$ that do not change sign such that for all $i=1,\dots,m$
\begin{eqnarray}\label{pointstability}
 \left\{ \begin{array}{lcl}
\hfill \div(y^{a_i} \nabla \phi_i)&=& 0   \ \ \text{in}\ \ \mathbf{R}_+^{n+1},\\   
\hfill -\lim_{y\to0}y^{a_i} \partial_{y} \phi_i&=& d_{s_i} \sum_{j=1}^m \partial_{v_i v_j}H(\mathbf v)\phi_j   \ \ \text{in}\ \ \partial\mathbf{R}_+^{n+1},
\end{array}\right.
  \end{eqnarray}
and $ \partial_{v_iv_j}H(\mathbf v) \phi_j \phi_i \ge 0$ for $1\leq i<j\leq m$.
 \item  A solution $\mathbf u=(u_i)_i$ of the system (\ref{main}) in $\mathbf{R}^n$ is said to be
{\it pointwise-stable}, if there exist $(\phi_i)_{i=1}^m$ in $C^\infty(\mathbf{R}^n)$ that do not change sign such that for all $i=1,\dots,m$
\begin{eqnarray}\label{pointstabilityu}
 (-\Delta)^{s_i}\phi_i=\sum_{j=1}^m \partial_{u_i u_j}H(\mathbf u)\phi_j   \ \ \text{in}\,\mathbf{R}^{n},
  \end{eqnarray}
and $ \partial_{u_iu_j}H(\mathbf u) \phi_j \phi_i \ge 0$ for $1\leq i<j\leq m$.
\end{enumerate}

 \end{dfn}
Indeed, consider a bounded pointwise-stable solution $\mathbf u$ of \eqref{main}. Then the extension $\mathbf v $ is constructed by convolving each component by the Poisson kernel associated to the operator $\mbox{div}(y^{1-2s_i}\nabla)$. Therefore, since the Poisson kernel is non negative, the function $\bf v$ will be bounded and pointwise-stable. 

\begin{dfn}
We say that $\mathbf v=(v_i)_i$ is a layer solution of (\ref{emain}) if $\mathbf v=(v_i)_i$ is a bounded solution of (\ref{emain}) such that for each $i$ the directional derivative $\partial_{x_n} v_i$ does not change sign and $\lim_{x_n\to \infty} \mathbf v=\mathbf\alpha$ and $\lim_{x_n\to -\infty} \mathbf v=\mathbf\beta$.  
\end{dfn}

\begin{notation} We fix the following notations throughout the paper. 
\begin{enumerate}
\item $s_*=\min_{i=1}^{m}\{s_i\}$,  $s^*=\max_{i=1}^{m}\{s_i\}$  and $0<s_*\le s^*<1$. 
\item $B_R^+=\{X=(\mathbf x,y)\in\mathbf R^{n+1}, |\mathbf X|<R\}$, $\partial^+ B_R^+=\partial B^+_R\cap \{y>0\}$ and $\Gamma^0_R=\partial B^+_R\cap \{y=0\}$. 
\item $C_R=B_R\times (0,R) \subset \mathbf R^{n+1}$. 
\end{enumerate} 
\end{notation}

\subsection{Main results} 
We can now state our main results. The following theorem is a symmetry result for $H-$monotone solutions of (\ref{emain}). Note that similar result holds for pointwise-stable and stable solutions. 
\begin{thm}\label{thsymv}
Let $n \leq 3$ and $1/2 \le s_* <1$. Suppose that ${\mathbf v}$ is a bounded $H$-monotone solution of the orientable system \eqref{emain}. Then, there exist a constant $\mathbf \Gamma_i\in\mathbf{S}^{n-1}$ and $v^*_i: \mathbf{R^+}\times  \mathbf{R}^+\to  \mathbf{R}$ such that
$$ v_i(\mathbf x,y)=v^*_i(\mathbf  \Gamma_i\cdot\mathbf x,y)$$
 for all $\mathbf x\in\mathbf{R}^{n}$ and $i=1,\cdots,m$.
\end{thm}

\begin{rem}
The case of $0 < s_* <1/2$ seems a challenging problem, even in the case of scalar equation, and will remain as an open problem.   
\end{rem}

The second main result is a Hamiltonian identity.

\begin{thm}\label{hamilton}
\begin{enumerate}
\item Let  $n=1$ and ${\bf v}$ be a layer solution of \eqref{emain}. Assume that $s_i=s \in (0,1)$ for every $i=1,...,m$. Then the following identity holds for any $x\in \mathbf R$
\begin{equation}\label{hamiltonian}
 \sum_{i=1}^{m} \int_0^\infty y^{1-2s} \left[  (\partial_x v_i)^2 - (\partial_y v_i)^2 \right] dy=2d_s \left[  H(\mathbf v(x,0)) - H(\mathbf \alpha) \right]. 
\end{equation}
\item Let $\bf v$ be a bounded radial solution of \eqref{emain}, i.e. $\mathbf{v} =\mathbf{v}(|\mathbf x|,y)$. Assume furthermore that $s_i=s$ for any $i$. Then, the following function is nonincreasing in $r$,
  \begin{equation}\label{radialmono}
   \sum_{i=1}^{m} \int_0^\infty y^{1-2s} \left[  (\partial_r v_i)^2 - (\partial_y v_i)^2 \right] dy - 2 d_s H(\mathbf v(r,0)).
  \end{equation}
\end{enumerate}
\end{thm}

The following theorem is a structural property on the nonlinearity to admit radial solutions. 
\begin{thm}\label{radialProp}
Assume that $n>1$ and $\mathbf v=(v_i)_i$ is a radial solution of (\ref{emain}) such that $\lim_{|\mathbf x|\to \infty} \mathbf v(|\mathbf x|,0)=0$. Then $\partial_{v_i} H(0)=0$ and $H(\mathbf v(0,0)) <H(0)$. Moreover, if $v_i(|\mathbf x|,y)$ is decreasing in $|\mathbf x|$ for all $1\le i\le m$ and $y>0$, then $$\sum_{i,j=1}^m \partial_{v_i v_j} H(0) \le 0.$$
\end{thm}
Finally, we present a Liouville theorem for nonlocal gradient systems with non sign changing nonlinearities in lower dimensions.     
 
\begin{thm}\label{lioupositive}
Suppose that $\mathbf v=(v_i)_i$ is a bounded pointwise-stable solution of (\ref{emain}) where $\nabla H \ge0$. Then, each $u_i$ is constant provided $n\le 2(1+s_*)$.  In particular when $0<s_*<1/2$ then $n\le 2$  and when $1/2\le s < 1$ then $n \le 3$.  
\end{thm}

\section{Preliminary results}
In this section we collect some useful results for our purposes.   Since the nonlinearity $H$ is in $C^{2,\gamma}(\mathbf R^m)$ and $\mathbf u$ is assumed to be bounded, the proof of the next lemma follows directly from the regularity result of Cabr\'{e}-Sire in \cite{CS1}. 
\begin{lemma}
Let $H$ be a $C^{2,\gamma}(\mathbf R^n)$ function with $\gamma >\max(0,1-2s_*)$. 
Then, any bounded solution of \eqref{main} is $C^{2,\beta}(\mathbf R^n)$ for some $0 <\beta < 1$ depending only on $s_*$ and $\gamma$.
\end{lemma}

Similarly, the next lemma is taken from \cite{CS1}. Under the assumption of boundedness of $\bf v$, it applies to the system (\ref{emain}).  

\begin{lemma}\label{regv}
Let $R>0$. Let $\varphi \in C^\sigma (\Gamma^0_{2R})$ for some $\sigma \in (0,1)$ and 
$u \in L^\infty(B^+_{2R}) \cap H^1(B^+_{2R},y^{1-2s})$ be a weak solution of
\begin{equation*}
\label{problemBR}
\begin{cases}
\mbox{div}(y^{1-2s}\nabla u)=0&\text{ in } B^+_{2R}\subset\mathbf R^{n+1}_+\\ 
-\lim_{y \to 0} y^{1-2s_i}\frac{\partial u}{\partial y}
=\varphi&\text{ on } \Gamma^0_{2R}.
\end{cases}
\end{equation*}
Then, there exists  $\beta \in (0,1)$ depending only on $n$, $s_i$, and $\sigma$, 
such that $u \in C^{0,\beta}(\overline{B_R^+})$ and 
$y^{1-2s_i} u_y \in C^{0,\beta}(\overline{B_R^+})$.  Furthermore, there exist constants $C^1_R$ and $C^2_R$ depending only on $n$, $a$, 
$R$, $\|u\|_{L^\infty(B_{2R}^+)}$ and also on $\|\varphi \|_{L^\infty(\Gamma^0_{2R})}$ (for $C^1_R$) 
and $\|\varphi \|_{C^\sigma(\Gamma^0_{2R})}$ (for $C^2_R)$ , such that 
$$\|u\|_{C^{0,\beta}(\overline{B_R^+})} \leq C^1_R \ \ 
\text{and} \ \ \|y^{1-2s_i} u_y\|_{C^{0,\beta}(\overline{B_R^+})} \leq C^2_R.$$
\end{lemma}

From the previous lemmata, one deduces the following gradient estimates. 

\begin{lemma}\label{asymp}
Let $\bf v$ be a bounded solution of \eqref{emain}. Then for every $i=1,\cdots,m$ one has
 \begin{eqnarray*}
|\nabla_{\mathbf x} v_i| &\le& \frac{C}{1+y} \ \ \text{for} \ \ (\mathbf x,y)\in\mathbf{R}^{n+1}_+ \\
|\partial_y v_i| &\le& \frac{C}{y}  \ \ \text{for} \ \ \mathbf x\in\mathbf{R}^{n}, y>1 \\
| y^{a_i} \partial_y v_i^t| &\le& C  \ \ \text{for} \ \  \mathbf x\in\mathbf{R}^{n}, 0<y<1 
  \end{eqnarray*}
\end{lemma}

The next lemma is a consequence of the previous regularity lemmata (as in \cite{CS1}). 
\begin{lemma}\label{decayv}
Suppose that the bounded vector function $\mathbf v$ is  either a layer solution of (\ref{emain}) in dimension $n=1$ or a radial solution of  (\ref{emain}) such that $\lim_{|\mathbf x|\to\infty} v(|\mathbf x|,0)$ exists then for each $1\le i\le m$, we have $$ \lim_{|\mathbf x|\to\infty} \int_0^\infty y^{a_i} |\nabla v_i|^2 dy=0.$$
\end{lemma}

As already mentioned, we investigate bounded solutions of \eqref{emain}. The following result is a consequence of the pointwise-stability. For the case of local elliptic systems, see \cite{FG}. 

\begin{lemma}\label{stabilitylem} 
Suppose that $\mathbf v=(v_i)_i$ is a pointwise-stable solution of system (\ref{emain}). Then the following stability inequality holds for orientable systems: 
\begin{equation}\label{stabilityin}
\sum_{i,j=1}^{m} \int_{\partial\mathbf{R}^{n+1}_+} \sqrt{d_{s_i} d_{s_j}}  \partial_{v_i v_j} H(\mathbf v)  \zeta_i \zeta_j  d\mathbf x  \le \sum_{i=1}^{m} \int_{\mathbf{R}^{n+1}_+}  y^{a_i} |\nabla \zeta_i|^2  d\mathbf x dy 
\end{equation}
for any $\zeta_i\in C_c^1(\mathbf R_+^{n+1})$.  
\end{lemma}

\noindent  \textbf{Proof:} Since  $\mathbf v=(v_i)_i$ is a pointwise-stable solution of system (\ref{emain}), there exists $\mathbf \phi=(\phi_i)_i$ that satisfies (\ref{pointstability}).  Multiple each equation by $\frac{\zeta_i^2}{\phi_i}$ and do integration by parts to get 
$$ - \int_{\mathbf R_+^{n+1}} y^{a_i} \nabla\phi_i\cdot \nabla \left(\frac{\zeta_i^2}{\phi_i} \right)+ \int_{\partial\mathbf R_+^{n+1}} y^{a_i} \nabla\phi_i\cdot \nu \left(\frac{\zeta_i^2}{\phi_i} \right)=0.$$
Using the boundary term of equation (\ref{pointstability}) we get 
 \begin{eqnarray}\label{}
- 2\int_{\mathbf R_+^{n+1}} y^{a_i} \nabla\phi_i  \cdot \nabla\zeta_i \frac{\zeta_i}{\phi_i} + \int_{\mathbf R_+^{n+1}} y^{a_i} |\nabla\phi_i|^2   \frac{\zeta_i^2}{\phi_i^2} + \int_{\partial\mathbf R_+^{n+1}} d_{s_i} \sum_{j=1}^m \partial_{v_i v_j}H(\mathbf v)\phi_j \frac{\zeta_i^2}{\phi_i}=0.
  \end{eqnarray}
Applying the Young's inequality, for any $i$, we obtain
 \begin{eqnarray}\label{}
\int_{\partial\mathbf R_+^{n+1}} d_{s_i} \sum_{j=1}^m \partial_{v_i v_j}H(\mathbf v) \frac{\phi_j}{\phi_i} \zeta_i^2&=&  \int_{\mathbf R_+^{n+1}} y^{a_i}  \left(   -  |\nabla\phi_i|^2   \frac{\zeta_i^2}{\phi_i^2} +2   \nabla\phi_i  \cdot \nabla\zeta_i \frac{\zeta_i}{\phi_i}  \right)
\\&\le&  \int_{\mathbf R_+^{n+1}} y^{a_i} |\nabla \zeta_i|^2
  \end{eqnarray}
Taking the sum on both sides for $i=1,\cdots,m$, the integrand in the left-hand side of the above inequality is 
 \begin{eqnarray*}
  \sum_{i,j} d_{s_i}  \partial_{v_i v_j}H(\mathbf v)  \frac{\phi_j}{\phi_i}\zeta_i^2 &=& \sum_{i} d_{s_i} \partial_{v_i v_j}H(\mathbf v)  \zeta_i^2 + \sum_{i\neq j} d_{s_i}  \partial_{v_i v_j}H(\mathbf v)   \frac{\phi_j}{\phi_i}\zeta_i^2
   \\&=& \sum_{i} d_{s_i} \partial_{v_i v_j}H(\mathbf v)  \zeta_i^2 + \sum_{i < j} d_{s_i} \partial_{v_i v_j}H(\mathbf v)   \frac{\phi_j}{\phi_i}\zeta_i^2 + \sum_{i> j} d_{s_i} \partial_{v_i v_j}H(\mathbf v)  \frac{\phi_j}{\phi_i}\zeta_i^2
    \\&= & \sum_{i} d_{s_i} \partial_{v_i v_j}H(\mathbf v)  \zeta_i^2 + \sum_{i < j} d_{s_i} \partial_{v_i v_j}H(\mathbf v) \frac{\phi_j}{\phi_i}\zeta_i^2 + \sum_{i< j} d_{s_j} \partial_{v_i v_j}H(\mathbf v)  \frac{\phi_i}{\phi_j}\zeta_j^2 
    \\ &=& \sum_{i} d_{s_i} \partial_{v_i v_j}H(\mathbf v)  \zeta_i^2 + \sum_{i < j} \partial_{v_i v_j}H(\mathbf v)  (\phi_i \phi_j)^{-1} \left( d_{s_i} \phi_j^2\zeta_i^2 +d_{s_j} \phi_i^2\zeta_j^2\right) 
    \\&\ge &  \sum_{i}  d_{s_i} \partial_{v_i v_j}H(\mathbf v)  \zeta_i^2 + 2\sum_{i < j}   \sqrt{d_{s_i}d_{s_j}} \partial_{v_i v_j}H(\mathbf v)  \zeta_i\zeta_j \ \ \text{since $\partial_{v_i v_j}H(\mathbf v)   (\phi_i \phi_j)^{-1} \ge 0$} \\&=&\sum_{i,j}   \sqrt{d_{s_i}d_{s_j}} \partial_{v_i v_j}H(\mathbf v) \zeta_i\zeta_j. 
  \end{eqnarray*} 
This finishes the proof. 

\hfill $ \Box$

\begin{dfn} A solution $\mathbf v$ of (\ref{emain}) that satisfies the stability inequality (\ref{stabilityin}) is called stable solution. 
\end{dfn}
\section{Symmetry of solutions and energy estimates}
In this section we prove a proof for Theorem \ref{thsymv}. Let us divide this section into two parts. A first one is devoted to the proof in dimension $n=2$. The result is known by \cite{serena} in the case of two equations. We generalize it to any number of equations, though the method is similar. The second part is devoted to the case of dimension $n=3$. 

\subsection{The case of dimension $n=2$}

We first derive a geometric Poincar\'e formula, which is valid in {\em any} dimension $n\ge 1$ and any number of equations $m\ge 1$. This is a generalization of the  geometric Poincar\'e formula given in \cite{sz} that is used by Farina-Sciunzi-Valdinoci \cite{FSV}, for the first time in the scalar case,  and by Fazly-Ghoussoub \cite{FG}, in the case of systems,  to prove certain De Giorgi type results. 

\begin{lemma}
 Assume that  $m,n\ge 1$ and $\mathbf v=(v_i)_i$ is a poinwise-stable solution of (\ref{emain}).  Then, for any $\eta=(\eta_k)_{k=1}^m \in C_c^1(\mathbf R_+^{n+1})$, the following inequality holds;
\begin{eqnarray}\label{poincare}
\nonumber
\sum_{i=1}^m \frac{1}{d_{s_i}}  \int_{\mathbf R_+^{n+1}} y^{a_i}  |\nabla_{\mathbf x}  v_i|^2   |\nabla \eta_i|^2 &\ge& \sum_{i=1}^m  \frac{1}{d_{s_i}}   \int_{    \{|\nabla_{\mathbf x}  v_i|\neq 0 \} \cap \mathbf R_+^{n+1} }   \left(   |\nabla v_i|^2 \mathcal{A}_i^2 + | \nabla_{T_i} |\nabla_{\mathbf x}  v_i| |^2  \right)\eta_i^2\\&&+\sum_{i\neq j} \int_{\partial\mathbf R_+^{n+1}}   \left(   |\nabla_{\mathbf x}  v_i|  |\nabla_{\mathbf x}  v_j| \eta_i \eta_j   - \nabla_{\mathbf x}  v_i \cdot   \nabla_{\mathbf x}  v_j \eta_i^2 \right)\partial_{v_i v_j}H(\mathbf v) ,
  \end{eqnarray} 
 where $\nabla_{T_i}$ stands for the tangential gradient along a given level set of $v_i$ and 
$\mathcal{A}_i^2$ for the sum of the squares of the principal curvatures of such a level set.
\end{lemma}
\noindent\textbf{Proof:} Since  $\mathbf v=(v_i)_i$ is a poinwise-stable solution of (\ref{emain}), from Lemma \ref{stabilitylem} the stability inequality  (\ref{stabilityin}) hold. Test the stability inequality with $\zeta_i:=|\nabla_{\mathbf x} v_i|\eta_i$ where each $\eta_i\in C^1_c(\mathbf R_+^{n+1})$, to get 
 \begin{eqnarray}\label{stablepoin}
I&:=& \sum_{i,j=1}^{m} \int_{\partial\mathbf{R}^{n+1}_+} \sqrt{d_{s_i} d_{s_j}}  \partial_{v_i v_j} H(\mathbf v)  |\nabla_{\mathbf x} v_i|  |\nabla_{\mathbf x} v_j|\eta_i \eta_j  d\mathbf x \\&\le& \sum_{i=1}^{m} \int_{\mathbf{R}^{n+1}_+}  y^{a_i} |\nabla \left(   |\nabla_{\mathbf x} v_i|\eta_i   \right) |^2  d\mathbf x dy =:J.
  \end{eqnarray} 
Simplifying $I $ we see that
 \begin{eqnarray}\label{Ipoin}
I:= \sum_{i=1}^{m} \int_{\partial\mathbf{R}^{n+1}_+} d_{s_i}   \partial_{v_i v_i} H(\mathbf v)  |\nabla_{\mathbf x} v_i|^2 \eta_i^2   d\mathbf x + 
\sum_{i\neq j} \int_{\partial\mathbf{R}^{n+1}_+} \sqrt{d_{s_i} d_{s_j}}  \partial_{v_i v_j} H(\mathbf v)  |\nabla_{\mathbf x} v_i|  |\nabla_{\mathbf x} v_j|   \eta_i \eta_j  d\mathbf x 
  \end{eqnarray} 
For the other term $J$, we have
\begin{eqnarray}\label{Jpoin}
J&=&  \sum_{i=1}^{m} \int_{\mathbf{R}^{n+1}_+}  y^{a_i} \left[   |\nabla_{\mathbf x} v_i|^2 |\nabla \eta_i|^2 +    |\nabla|\nabla_{\mathbf x} v_i||^2  \eta_i^2  +\frac{1}{2}   \nabla |\nabla_{\mathbf x} v_i|^2  \cdot \nabla \eta_i^2  \right] d\mathbf x dy. 
  \end{eqnarray} 
Differentiate the $i^{th}$ equation of (\ref{emain}) with respect to $x_k$ for each $i=1,2,...,m$ and $k=1,2,...,n$ and multiply with $\partial_{x_k} v_i \eta_i^2$ to get 
$$ \div(y^a_i   \nabla \partial_{x_k} v_i  ) \partial_{x_k} v_i \eta_i^2=0. $$
Integrate by parts to obtain
 \begin{eqnarray}\label{partsxk}
\int_{\mathbf{R}^{n+1}_+}  y^{a_i}  |\nabla \partial_{x_k} v_i|^2 \eta_i^2 + \frac{1}{2} \int_{\mathbf{R}^{n+1}_+}  y^{a_i} \nabla |\partial_{x_k} v_i|^2 \cdot \nabla \eta_i^2 = \int_{\partial\mathbf{R}^{n+1}_+} \lim_{y\to 0} y^{a_i}  (-\partial _y \partial_{x_k} v_i) \partial_{x_k} v_i \eta_i^2
  \end{eqnarray} 
Differentiating the boundary term in of (\ref{emain}) with respect to $x_k$ we get 
$$ \lim_{y\to 0} y^{a_i}  (-\partial _y \partial_{x_k} v_i) \partial_{x_k} v_i \eta_i^2 = d_{s_i} \sum_{j=1}^m \partial_{v_i v_j} H(\mathbf v) \partial_{x_k} v_j \partial_{x_k} v_i   \eta_i^2.$$
From this and (\ref{partsxk}) and taking sum on $i$ we have
 \begin{eqnarray}\label{partsxkk}
\nonumber&& \sum_{i=1}^{m} \int_{\mathbf{R}^{n+1}_+}  y^{a_i}  |\nabla \partial_{x_k} v_i|^2 \eta_i^2 +  \sum_{i=1}^{m} \frac{1}{2} \int_{\mathbf{R}^{n+1}_+}  y^{a_i} \nabla |\partial_{x_k} v_i|^2 \cdot \nabla \eta_i^2  \\&=&  \sum_{i=1}^{m} d_{s_i} \int_{\partial\mathbf{R}^{n+1}_+}  \partial_{v_i v_i} H(\mathbf v) |\partial_{x_k} v_i|^2   \eta_i^2   +  \sum_{j=1}^m  d_{s_i}  \int_{\partial\mathbf{R}^{n+1}_+}  \partial_{v_i v_j} H(\mathbf v) \partial_{x_k} v_j \partial_{x_k} v_i   \eta_i^2  
  \end{eqnarray} 
Taking sum on $k=1,\cdots, n$ we have 
 \begin{eqnarray}\label{partsumk}
 \sum_{i=1}^{m} d_{s_i} \int_{\partial\mathbf{R}^{n+1}_+}  \partial_{v_i v_i} H(\mathbf v) |\nabla_{\mathbf x}  v_i|^2   \eta_i^2 &=&  \sum_{i=1}^{m} \int_{\mathbf{R}^{n+1}_+}  y^{a_i}   \sum_{k=1}^{n}   |\nabla \partial_{x_k} v_i|^2 \eta_i^2
  \\ && \nonumber + \frac{1}{2}  \sum_{i=1}^{m} \int_{\mathbf{R}^{n+1}_+}  y^{a_i} \nabla |\nabla_{\mathbf x} v_i|^2 \cdot \nabla \eta_i^2  \\ &&  \nonumber -    \sum_{j \neq i}^m  d_{s_i}  \int_{\partial\mathbf{R}^{n+1}_+}  \partial_{v_i v_j} H(\mathbf v)   \nabla_{\mathbf x} v_i \cdot \nabla_{\mathbf x} v_j  \eta_i^2 
  \end{eqnarray} 
Substitute (\ref{partsumk}) in $I\le J$ where $I$ is given in (\ref{Ipoin}) and $J$ is given in (\ref{Jpoin}) then the term $$\frac{1}{2}  \sum_{i=1}^{m} \int_{\mathbf{R}^{n+1}_+}  y^{a_i} \nabla |\nabla_{\mathbf x} v_i|^2 \cdot \nabla \eta_i^2$$ will be cancel out and we end up with 
  \begin{eqnarray}\label{partsumk}
&&\sum_{i=1}^{m} \int_{\mathbf{R}^{n+1}_+ \cap \{ | \nabla v_i | \neq 0 \} }  y^{a_i}  \left( \sum_{k=1}^{n}   |\nabla \partial_{x_k} v_i|^2   -     |\nabla|\nabla_{\mathbf x} v_i||^2      \right) \eta_i^2
\\&& +  \sum_{i\neq j} \int_{\partial\mathbf{R}^{n+1}_+}     \partial_{v_i v_j} H(\mathbf v) \left[  \sqrt{d_{s_i} d_{s_j}}   |\nabla_{\mathbf x} v_i|  |\nabla_{\mathbf x} v_j|   \eta_i \eta_j -  d_{s_i}  \partial_{v_i v_j} H(\mathbf v)   \nabla_{\mathbf x} v_i \cdot \nabla_{\mathbf x} v_j \eta_i^2  \right]
\\ &\le&  \sum_{i=1}^{m} \int_{\mathbf{R}^{n+1}_+ }  y^{a_i} |\nabla_{\mathbf x} v_i|^2 |\nabla \eta_i|^2
 \end{eqnarray} 
According to formula (2.1) given in  \cite{sz}, the following geometric identity between the tangential gradients and curvatures holds. For any $w \in C^2(\Omega)$
 \begin{eqnarray}\label{identity}
  \sum_{k=1}^{n} |\nabla \partial_k w|^2-|\nabla|\nabla w||^2=
\left\{
                      \begin{array}{ll}
                       |\nabla w|^2 (\sum_{l=1}^{n-1} \mathcal{\kappa}_l^2) +|\nabla_T|\nabla w||^2 & \hbox{for $x\in\{|\nabla w|>0\cap \Omega \}$,} \\
                       0 & \hbox{for $x\in\{|\nabla w|=0\cap \Omega \}$,}
                                                                       \end{array}
                    \right.   \end{eqnarray} 
 where $ \mathcal{\kappa}_l$ are the principal curvatures of the level set of $w$ at $\mathbf x$ and $\nabla_T$ denotes the orthogonal projection of the gradient along this level set . This completes the proof.
     
\hfill $ \Box$

The proof of next lemma is straightforward and it is omitted.

\begin{lemma} Suppose that $\mathbf \phi=(\phi_i)_i$ and  $\mathbf \psi=(\psi_i)_i$ are solutions of the linear equation
\begin{eqnarray}\label{pointstability1}
 \left\{ \begin{array}{lcl}
\hfill \div(y^{a_i} \nabla w_i)&=& 0   \ \ \text{in}\ \ \mathbf{R}_+^{n+1},\\   
\hfill -\lim_{y\to0}y^{a_i} \partial_{y} w_i&=& d_{s_i} \sum_{j=1}^m \partial_{v_i v_j}H(\mathbf v)w_j   \ \ \text{in}\ \ \partial\mathbf{R}_+^{n+1}.
\end{array}\right.
  \end{eqnarray}
  Assume that $\phi_i$ does not change sign for each $i$. Then,  $\mathbf \sigma=(\sigma_i)_i$ where the quotient $\sigma_i:=\frac{\psi_i}{\phi_i}$ satisfies  
  \begin{eqnarray}\label{sigma}
 \left\{ \begin{array}{lcl}
\hfill \div(y^{a_i} \phi_i^2 \nabla \sigma_i)&=& 0   \ \ \text{in}\ \ \mathbf{R}_+^{n+1},\\   
\hfill -\lim_{y\to0}y^{a_i} \phi_i^2 \partial_{y} \sigma_i&=& d_{s_i} \sum_{j=1}^m \partial_{v_i v_j}H(\mathbf v) \phi_i \phi_j   (\sigma_j-\sigma_i) \ \ \text{in}\ \ \partial\mathbf{R}_+^{n+1}.
\end{array}\right.
  \end{eqnarray}
\end{lemma}

The next theorem is crucial to get the symmetry result in dimension $n=2$. The proof is very similar to the case of $m=2$ given in \cite{serena} and we omit it here. 
 \begin{thm}\label{symn2}
Suppose that $\mathbf v=(v_i)_i$ is a bounded stable solution of (\ref{emain}) such that 
$$ \sum_{i=1}^{m}  \int_{B_R^+}  y^{a_i} |\nabla v_i|^2 \le C R^2 $$
where $C$ is independent of $R$. Then, there exists a constant $\mathbf\Lambda_i\in\mathbf{S}^{n-1}$ and $v^*_i: \mathbf{R^+}\times  \mathbf{R}\to  \mathbf{R}$ such that
$$ v_i(\mathbf x,y)=v^*_i(\mathbf\Lambda_i\cdot\mathbf x,y)$$
 for all $(\mathbf x,y)\in\mathbf{R}^{n+1}$ and $i=1,\cdots,m$.
\end{thm}

\subsection{The case of dimension $n=3$}

\subsubsection{A linear Liouville theorem}\label{sec1}
In this part, we provide a linear Liouville theorem for an elliptic system in half-space. On the entire space $\mathbf R^n$ this theorem is given in \cite{FG}. For the scalar case, $m=1$,  this type of Liouville theorem was noted by Berestycki, Caffarelli and Nirenberg in \cite{BCN1} and used by Ghoussoub-Gui \cite{GG1} and later by Ambrosio and Cabr\'{e} \cite{AC} to prove the De Giorgi conjecture in dimensions two and three. Then this theorem was improved by Moschini \cite{Mos}. For nonlocal equations,  this type linear Liouville theorem is given by Cabr\'{e}-Sire \cite{CS}.     Consider the set of functions with a limited growth as infinity as $$\mathcal F=\left\{F:\mathbf R^+\to\mathbf R^+, F \  \text{is nondecreasing and} \ \int_{2}^{\infty} \frac{1}{rF(r)}=\infty\right\}.$$
In particular, $F(r)=\log r$ belongs to this class and $F(r)=r$ does not belong to $\mathcal F$. As far as we know,  this class of functions was considered by Karp \cite{Karp1,Karp2} for the first time. Here is the Liouville theorem. 

\begin{thm}\label{liouville} Assume that for each $i=1,\cdots,m$,  $\phi_i \in L^{\infty}_{loc} (\overline{\mathbf{R}_+^{n+1}}) $ is a positive function and $\sigma_i \in H^1_{loc}(\overline{\mathbf{R}_+^{n+1}},y^{a_i}) $ satisfy
\begin{equation}\label{liouassum}
 \limsup_{R\to\infty} \frac{1}{R^2F(R)} \int_{C_R}  \sum_{i=1}^{m}  y^{a_i}  \phi_i^2\sigma_i^2 <\infty.
 \end{equation}
 If $(\sigma_i)_{i=1}^m$ are solutions of 
 \begin{eqnarray}\label{div}
 \left\{ \begin{array}{lcl}
\hfill  -\sigma_i \div(y^{a_i}\phi_i^2\nabla \sigma_i)&\le& 0   \ \ \text{in}\ \ \mathbf{R}_+^{n+1},\\   
\hfill -\lim_{y\to0}y^{a_i} \phi_i \sigma_i \partial_{y} \sigma_i &\le&  \sum_{j=1}^{n} h_{ij} f(\sigma_j-\sigma_i)\sigma_i  \ \ \text{in}\ \ \partial\mathbf{R}_+^{n+1},
\end{array}\right.
  \end{eqnarray}
where $0\le h_{i,j}\in L_{loc}^1(\mathbf{R}^n)$, $h_{i,j}=h_{j,i}$ and $f\in L_{loc}^1(\mathbf{R})$ is an odd function such that $f(t)\ge 0$ for $t\in\mathbf R^+$.  Then, for all  $i=1,\cdots,m$, the functions $\sigma_i$ are constant.
\end{thm}
\noindent\textbf{Proof:} Note that using the equation (\ref{div}) we obtain
 \begin{eqnarray*}
\div(y^{a_i} \phi_i^2 \sigma_i \nabla \sigma_i) &= & \div(y^{a_i} \phi_i^2  \nabla \sigma_i) \sigma_i + y^{a_i} \phi_i^2  |\nabla \sigma_i |^2 \\&\ge&  y^{a_i} \phi_i^2  |\nabla \sigma_i |^2. 
  \end{eqnarray*}
Set $C_R=B_R \times (0,R)$. Integrating the above on $C_R$ and using the boundary term in (\ref{div}) we get 
\begin{eqnarray*}
I:= \int_{C_R} \sum_{i=1}^{m} y^{a_i} \phi_i^2  |\nabla \sigma_i |^2 &\le& 
 \int_{C_R} \sum_{i=1}^{m} \div(y^{a_i} \phi_i^2 \sigma_i \nabla \sigma_i) 
 \\&=& \int_{\partial^+ C_R} \sum_{i=1}^{m} y^{a_i} \phi_i^2 \sigma_i \partial_\nu \sigma_i 
 +\int_{ B_R\times \{y=0\}} \sum_{i=1}^{m} y^{a_i} \phi_i^2 \sigma_i (-\partial_y \sigma_i) 
 \\&\le&   \int_{\partial^+ C_R} \sum_{i=1}^{m} y^{a_i} \phi_i^2 \sigma_i \partial_\nu \sigma_i +\int_{ B_R\times \{y=0\}} \sum_{i,j=1}^{m} h_{ij} f(\sigma_j-\sigma_i) \sigma_i 
 \\&=:& J(R)+K(R). 
  \end{eqnarray*}
Note that $K\le 0$ since 
\begin{eqnarray*}
 \sum_{i,j} h_{ij}(x) \sigma_i f(\sigma_j-\sigma_i)&=& \sum_{i< j}  h_{i,j} \sigma_i f(\sigma_j-\sigma_i) + \sum_{i> j} h_{i,j}  \sigma_i  f(\sigma_j-\sigma_i)  \\&=&\sum_{i< j} h_{i,j}  \sigma_i f(\sigma_j-\sigma_i)  + \sum_{i< j}h_{i,j} \sigma_j f(\sigma_i-\sigma_j)   \ \ \text{since} \ \ h_{ij}=h_{ji}\\&=&-\sum_{i< j} h_{i,j} (\sigma_j-\sigma_i)f(\sigma_j-\sigma_i)  \ \ \text{because $f$ is odd}.
  \end{eqnarray*}
Note also that 
\begin{eqnarray*}
J (R)&\le& \int_{\partial^+ C_R}   \left( \sum_{i=1}^{m} y^{a_i} \phi_i^2 \sigma_i^2  \right)^{1/2}  \left( \sum_{i=1}^{m} y^{a_i} \phi_i^2 |\nabla \sigma_i|^2  \right)^{1/2}
\\&\le&  \left(  \int_{\partial^+ C_R} \sum_{i=1}^{m} y^{a_i} \phi_i^2 \sigma_i^2 \right)^{1/2}  \left(  \int_{\partial^+ C_R} \sum_{i=1}^{m} y^{a_i} \phi_i^2 |\nabla \sigma_i|^2 \right)^{1/2}
\\&=& (I'(R))^{1/2}  \left(  \int_{\partial^+ C_R} \sum_{i=1}^{m} y^{a_i} \phi_i^2 \sigma_i^2 \right)^{1/2}  
  \end{eqnarray*}
From this we obtain  $$ I(R)\le J(R) \le  (I'(R))^{1/2}  \left(  \int_{\partial^+ C_R} \sum_{i=1}^{m} y^{a_i} \phi_i^2 \sigma_i^2 \right)^{1/2} $$
Therefore, 
\begin{equation}
\left( \int_{\partial^+ C_R} \sum_{i=1}^{m} y^{a_i} \phi_i^2 \sigma_i^2\right)^{-1} \le \frac{I'(R)}{I^2(R)}
\end{equation}
Suppose that at least one of the $\sigma_i$ is not constant. Then there exists $R_0>0$ such that $I(R)>0$ for every $R>R_0$. 
\begin{eqnarray*}
(r_2-r_1)^2 \left(  \int_{C_{r_2} \setminus C_{r_1}}  \sum_{i=1}^{m} y^{a_i} \phi_i^2 \sigma_i^2  \right)^{-1}
&=& (r_2-r_1)^2  \left(  \int_{r_1}^{r_2} dR \int_{\partial^+ C_R}  \sum_{i=1}^{m} y^{a_i} \phi_i^2 \sigma_i^2  \right)^{-1} \\& \le& \int_{r_1}^{r_2} dR \left( \int_{\partial^+ C_R} \sum_{i=1}^{m} y^{a_i} \phi_i^2 \sigma_i^2  \right)^{-1}  \\& \le& \frac{1}{I(r_1)} -  \frac{1}{I(r_2)}
  \end{eqnarray*}
Set $r_2=2^{j+1} r_*$ and $r_1=2^j r_*$ for some $r_*>R_0$ when $j=0,\cdots,M$ we have
$$ \sum_{j=1}^{m} \frac{2^{2j}}{   \int_{C_{2^{j+1} r_*}}   \sum_{i=1}^{m} y^{a_i} \phi_i^2 \sigma_i^2   } \le \frac{1}{I(r_*)}.$$
 From this and the assumption (\ref{liouassum}) we obtain
 $$  \sum_{j=1}^{m} \frac{1}{F(2^{j+1}r_*)} \le C  \sum_{j=1}^{m} \frac{2^{2j}}{   \int_{C_{2^{j+1} r_*}}   \sum_{i=1}^{m} y^{a_i} \phi_i^2 \sigma_i^2   } \le \frac{1}{I(r_*)}. $$
If we send $m\to\infty$, one can see that the series is divergent that is a contradiction. 

\hfill $ \Box$

\subsubsection{Energy estimate for layer solutions}

The energy functional is given by 
\begin{equation}\label{energy}
E_R(\mathbf v)= \sum_{i=1}^{m}\frac{1}{2d_{s_i}} \int_{C_R}  y^{a_i} |\nabla v_i|^2 -  \int_{B_R} H( \mathbf v) 
\end{equation}
where $C_R=B_R\times (0,R)$. then we derive energy estimates for monotone solutions of (\ref{main}).
\begin{thm}\label{energylayer}
 Suppose that $\mathbf v=(v_i)_i$ is a bounded $H$-monotone solution of (\ref{emain}) such that   $$\lim_{x_n\to\infty} \mathbf v(\mathbf x',x_n)=\mathbf \alpha$$ where $\mathbf \alpha=(\alpha_i)_i \in\mathbf{R}^n$ and $H(\mathbf \alpha)=0$. Then for various values of $s_*$, the following holds.
\begin{enumerate}
\item If $1/2<s_*<1$, then 
\begin{equation*}
E_R(\mathbf v) \le  C R^{n-1}.
\end{equation*}
\item For the case $s_*=1/2$, the energy bound is
\begin{equation*}
E_R(\mathbf v) \le  C R^{n-1}\log R. 
\end{equation*}
\item When $0<s_*<1/2$, then 
\begin{equation*}
E_R(\mathbf v) \le  C R^{n-2s_*}.
\end{equation*}
\end{enumerate} 
\end{thm}

 \noindent  \textbf{Proof:}  Define the shift function $v_i^t(\mathbf x,y):=v_i(\mathbf x',x_n+t,y)$ for $(\mathbf x',x_n,y)\in\mathbf R^{n+1}$ and $t\in\mathbf R$. For each $i$, we have 
 \begin{eqnarray}\label{vt}
 \left\{ \begin{array}{lcl}
\hfill \div(y^{a_i} \nabla v^t_i)&=& 0   \ \ \text{in}\ \ \mathbf{R}_+^{n+1},\\   
\hfill -\lim_{y\to 0}y^{a_i} \partial_{y} v^t_i&=& d_{s_i}\partial_{v_i}H(\mathbf v^t)   \ \ \text{in}\ \ \partial\mathbf{R}_+^{n+1},
\end{array}\right.
  \end{eqnarray}
where $|\mathbf v^t|\in {L^\infty}$ for every $t$. Furthermore $\mathbf v^t$ satisfies the pointwise bounds in Lemma \ref{asymp} and for all $(\mathbf x,y)\in\mathbf{R}^{n+1}_+$ we have $$ \lim_{t\to\infty} |v_i^t(\mathbf x,y)-\alpha_i|+|\nabla v_i^t(\mathbf x,y)|=0.$$
Therefore, $\lim_{t\to\infty} E_{R}(\mathbf v^t)=0$ where $\mathbf v^t=(v^t_1,\cdots,v^t_m)$. Differentiating the energy functional with respect to $t$ we get
 \begin{eqnarray*}
\partial_t E_R(\mathbf v^t)&=&  \sum_{i=1}^{m} \frac{1}{d_{s_i}}\int_{0}^{R} dy \int_{B_R}  y^{a_i} \nabla v^t_i\cdot\nabla(\partial_t v_i^t)-  \int_{B_R} \sum_{i=1}^{m}  H_{v_i}(\mathbf v^t) \partial_t v_i^t 
\\ &=&  -\sum_{i=1}^{m}\frac{1}{d_{s_i}}  \int_{0}^{R} dy \int_{B_R} \div( y^{a_i} \nabla v^t_i) \partial_t v_i^t + \sum_{i=1}^{m} \frac{1}{d_{s_i}} \int_{0}^{R} dy \int_{\partial B_R} y^{a_i} \nabla v^t_i\cdot\nu \partial_t v_i^t 
\\&&-  \sum_{i=1}^{m}  \int_{B_R} H_{v_i}(\mathbf v^t) \partial_t v_i^t 
  \end{eqnarray*}
Using (\ref{vt}) we have 
 \begin{eqnarray*}
 \partial_t E_R(\mathbf v^t)&=&\sum_{i=1}^{m} \frac{1}{d_{s_i}} \int_{0}^{R} dy \int_{\partial B_R} y^{a_i} \nabla v^t_i\cdot\nu \partial_t v_i^t - \sum_{i=1}^{m}  \int_{B_R} H_{\mathbf v_i}(v^t) \partial_t v_i^t \\&=&
 \sum_{i=1}^{m}\frac{1}{d_{s_i}} \left(  \int_{0}^{R} dy \int_{\partial B_R} y^{a_i} \nabla v^t_i\cdot\nu \partial_t v_i^t + \int_{B_R\times\{y=R\}}    y^{a_i} \partial_y v^t_i \partial_t v_i^t  \right)
   \end{eqnarray*}
Let us fix the following indices, $$ \partial_t v_\mu ^t >0> \partial_t v_\lambda ^t$$
when $\mu\in I$ and $\lambda\in J$ where $I\cup J=\{1,\cdots,m\}$. Therefore,
 \begin{eqnarray*}
 \partial_t E_R(\mathbf v^t)&=& \int_{\partial B_R}\int_{0}^{R} dy \left( \sum_{\mu\in I} \frac{1}{d_{s_\mu}} y^{a_\mu} \partial_\nu v^t_\mu  \partial_t v_\mu^t + \sum_{\lambda\in J} \frac{1}{d_{s_\lambda}} y^{a_\lambda} \partial_\nu v^t_\lambda  \partial_t v_\lambda^t  \right)
 \\&&+ \int_{B_R\times\{y=R\}} \left( \sum_{\mu\in I} \frac{1}{d_{s_\mu}} y^{a_\mu} \partial_y v^t_\mu  \partial_t v_\mu^t + \sum_{\lambda\in J}  \frac{1}{d_{s_\lambda}} y^{a_\lambda} \partial_y v^t_\lambda  \partial_t v_\lambda^t  \right)
   \end{eqnarray*}
Since $|\nabla v^t_i|\in L^\infty(\mathbf R^n)$ for each $i$, we have $-\frac{M}{1+y}\le \partial_\nu v^t_i \le \frac{M}{1+y}$ for $(x,y)\in\mathbf{R}^{n+1}_+$ and $-\frac{C}{y} \le \partial_y v_i^t \le \frac{C}{y}$  for  $x\in\mathbf{R}^{n}$ and $ y>1$. Using these we get the following. 
\begin{eqnarray*}
 \partial_t E_R(\mathbf v^t)&\ge&  \int_{\partial B_R}\int_{0}^{R} dy \left( \sum_{\mu\in I}  \frac{1}{d_{s_\mu}} y^{a_\mu} \left(- \frac{M}{1+y} \right)  \partial_t v_\mu^t + \sum_{\lambda\in J} \frac{1}{d_{s_\lambda}}  y^{a_\lambda} \left( \frac{M}{1+y} \right)  \partial_t v_\lambda^t  \right)
 \\&&+ \int_{B_R\times\{y=R\}} \left( \sum_{\mu\in I}  \frac{1}{d_{s_\mu}}  y^{a_\mu} \left(- \frac{C}{y} \right)  \partial_t v_\mu^t + \sum_{\lambda\in J} \frac{1}{d_{s_\lambda}}  y^{a_\lambda} \left(\frac{C}{y} \right)  \partial_t v_\lambda^t  \right)
   \end{eqnarray*}
Note that for every $T>0$ we have 
$$ E_R(\mathbf v)= E_R(\mathbf v^T)-\int_0^T \partial_t E_R(\mathbf v^t) dt.$$
Therefore, 
\begin{eqnarray*}
E_R(\mathbf v) &\le& E_R(\mathbf v^T) + M \int_{\partial B_R}\int_{0}^{R} dy \left( \sum_{\mu\in I} \frac{1}{d_{s_\mu}} \left(\frac{y^{a_\mu} }{1+y} \right)  \int_0^T \partial_t v_\mu^t dt - \sum_{\lambda\in J} \frac{1}{d_{s_\lambda}} \left( \frac{y^{a_\lambda} }{1+y} \right) \int_0^T \partial_t v_\lambda^t dt  \right)
 \\&&+ C\int_{B_R\times\{y=R\}} \left( \sum_{\mu\in I}\frac{1}{d_{s_\mu}} y^{a_\mu-1}   \int_0^T \partial_t v_\mu^t dt  - \sum_{\lambda\in J}  \frac{1}{d_{s_\lambda}} y^{a_\lambda-1}    \int_0^T \partial_t v_\lambda^t  dt \right).
   \end{eqnarray*}
Therefore, we can simplify the above as 
\begin{eqnarray*}
E_R(\mathbf v) &\le& E_R(\mathbf v^T) + C \int_{\partial B_R} \int_{0}^{R} \left( \sum_{\mu\in I}  \frac{y^{a_\mu} }{1+y}      (v_\mu^T - v_\mu) + \sum_{\lambda\in J}  \frac{y^{a_\lambda} }{1+y}   (v_\lambda-v_\lambda^T) \right) dy
 \\&&+ C\int_{B_R\times\{y=R\}} \left( \sum_{\mu\in I} y^{a_\mu-1}   (v_\mu^T -v_\mu) + \sum_{\lambda\in J} y^{a_\lambda-1}   (v_\lambda-v_\lambda^T) \right).
   \end{eqnarray*}
Note that $v_\mu^T \ge v_\mu$ and $v_\lambda \ge v_\lambda^T$ and $v_i\in L^{\infty}$. So, 
 \begin{eqnarray*}
E_R(v) \le  E_R(v^T) + C \int_{\partial B_R}  \int_0^R  \sum_{i=1}^{m} \frac{y^{a_i} }{1+y} dy + C  \int_{B_R\times\{y=R\}}  \sum_{i=1}^{m} y^{a_i-1} 
   \end{eqnarray*}
   Now sending $T\to\infty$ we get $E_R(v^T)=0$. Doing the integration we get 
   \begin{eqnarray*}
E_R(\mathbf v) &\le&  C  \sum_{i=1}^{m}\left(R^{n-1+a_i} \chi_{\{0<a_i<1\}} +  R^{n-1} \chi_{\{-1<a_i<0\}} +   R^{n-1}\log R \chi_{\{a_i=1/2\}}   +   R^{n+a_i-1} \right) \\
&=&C  \sum_{i=1}^{m}\left(R^{n-2s_i} \chi_{\{0<s_i<1/2\}} +  R^{n-1} \chi_{\{1/2<s_i<1\}} +   R^{n-1}\log R \chi_{\{s_i=1/2\}}   +   R^{n-2s_i} \right) 
   \end{eqnarray*}   
Considering the notation $s_*:=\min_{i=1}^{m} \{s_i\}$, we get the desired energy estimates. 

               \hfill $ \Box$

\noindent \textbf{Proof of Theorem \ref{thsymv}:}    Let again $\phi_i := \partial _n u_i$ and $\psi_i:=\nabla u_i\cdot\eta$ for any fixed $\eta=(\eta',0)\in \mathbf{R}^{n-1}\times\{0\}$ in such a way 
that $\sigma_i:=\frac{\psi_i}{\phi_i}$  is a solution of system (\ref{sigma}).  Set  $h_{i,j}(\mathbf x)=H_{u_iu_j}\phi_i(\mathbf x)\phi_j(\mathbf x)$ and $f$ to be the identity in Theorem \ref{liouville}. Note that $\mathbf u$ is a $H$-monotone solution, so $h_{i,j}\le 0$.  In dimension $n=2$,  assumption (\ref{liouassum}) holds and Theorem \ref{liouville} then yields that $\sigma_i$  is constant, which finishes the proof as argued before. 

 In dimension $n=3$, we shall follow methods developed by  Ambrosio-Cabr\'{e} \cite{AC} and Alberti-Ambrosio-Cabr\'{e} \cite{AAC}, in the case of local equations,  that is used later by Fazly-Ghoussoub  \cite{FG} in the case of local systems and by Cabr\'{e}-Cinti  \cite{cabre}  in the case of nonlocal equations.
 
 We first note that $\mathbf v$ is a bounded stable solution of (\ref{emain}) in $\mathbf{R}^4_+$. So,  the function $\bar {\mathbf v}(x_1, x_2,y):=\lim_{x_3\to \infty} \mathbf v(x_1, x_2, x_3,y)$ is also a bounded stable solution for (\ref{emain}) in $\mathbf{R}^3_+$.  It follows from Theorem \ref{symn2} that $\mathbf v$ is one dimensional and consequently the energy of $\bar {\mathbf v}$ in $C_R\subset \mathbf{R}^3_+ $ is bounded by a multiple of $R$, which yields that
 \begin{equation}
\label{boundEt} 
  E_R( \bar{\mathbf v}) \le C R^{2}, 
 \end{equation}
 where here $E_R(\mathbf v)=\sum_{i=1}^{m}\frac{1}{2d_{s_i}} \int_{C_R}  y^{a_i} |\nabla v_i|^2 d\mathbf x dy - \int_{B_R} \left( H( \mathbf v) -c_{\mathbf v} \right)d\mathbf x$ for $c_{\mathbf v}:=\sup  H(\mathbf v)$.    To finish the proof, we shall show that  
\begin{equation}
\label{ebound}
\sum_{i=1}^{m}\frac{1}{2d_{s_i}} \int_{C_R}  y^{a_i} |\nabla v_i|^2 d\mathbf x dy \le  C R^{2}  \ \chi_{\{s_*>1/2\}} + C  R^{2}\log R\  \chi_{\{s_*=1/2\}} .
 \end{equation}
Note that shifted function $\mathbf v^t(\mathbf x',y):=\mathbf v(\mathbf x',x_n+t,y)$ for $t\in\mathbf{R}$ and is also a bounded solution of (\ref{emain}), i.e.,
\begin{eqnarray}\label{emaint}
 \left\{ \begin{array}{lcl}
\hfill \div(y^{a_i} \nabla v^t_i)&=& 0   \ \ \text{in}\ \ \mathbf{R}_+^{n+1}=\left \{x \in \mathbf R^n, y>0 \right \},\\   
\hfill -\lim_{y\to0}y^{a_i} \partial_{y} v^t_i&=& d_{s_i} \partial_{v^t_i}H(\mathbf v^t)   \ \ \text{in}\ \ \partial\mathbf{R}_+^{n+1}.
\end{array}\right.
  \end{eqnarray}
Since  $v_i^t$ converges to $\bar v_i$ in $C^1_{loc}(\mathbf{R}^n)$ for all $i=1,\cdots,m$, we have $$ \lim_{t\to\infty} E_R(\mathbf v^t)= E_R(\bar{\mathbf v}).
$$
Note that for every $T>0$ we have 
$$ E_R(\mathbf v)= E_R(\mathbf v^T)-\int_0^T \partial_t E_R(\mathbf v^t) dt.$$
Similar to the proof of Theorem \ref{energylayer} we have 
 \begin{eqnarray*}
E_R(\mathbf v) \le  E_R(\mathbf v^T) + C \int_{\partial B_R}  \int_0^R  \sum_{i=1}^{m} \frac{y^{a_i} }{1+y} dy + C  \int_{B_R\times\{y=R\}}  \sum_{i=1}^{m} y^{a_i-1} 
   \end{eqnarray*}
  Sending $T\to\infty$ we get 
    \begin{eqnarray*}
E_R(\mathbf v) &\le&  E_R(\bar{\mathbf v}) + C \int_{\partial B_R}  \int_0^R  \sum_{i=1}^{m} \frac{y^{a_i} }{1+y} dy + C  \int_{B_R\times\{y=R\}}  \sum_{i=1}^{m} y^{a_i-1} \\
&\le & C R^2 + C \int_{\partial B_R}  \int_0^R  \sum_{i=1}^{m} \frac{y^{a_i} }{1+y} dy + C  \int_{B_R\times\{y=R\}}  \sum_{i=1}^{m} y^{a_i-1}
\\&\le& C R^2 +C  \sum_{i=1}^{m}\left(R^{n-2s_i} \chi_{\{0<s_i<1/2\}} +  R^{n-1} \chi_{\{1/2<s_i<1\}} +   R^{n-1}\log R \chi_{\{s_i=1/2\}}   +   R^{n-2s_i} \right) 
   \end{eqnarray*}  
Suppose that $1/2 < s_*<1$ and $n=3$ . Then for any $1\le i\le m$, we have $n-2s_i\le n-2s_* < n-1$ and therefore 
$$E_R(\mathbf v)  \le  C R^2 +C  R^{n-1} \log R \le C R^2 .$$
For the case $s_*=1/2$ and $n=3$, we have 
$$E_R(\mathbf v)  \le  C R^2 +C  R^{n-1} \log R \le C R^2 \log R.$$
This finishes the the proof of (\ref{ebound}).

               \hfill $ \Box$

\section{Optimality of energy estimates via a monotonicity formula}
Following ideas provided by Cabr\'{e}-Cinti \cite{cabre} for the scalar case and applying the Pohozaev identity, we prove the following monotonicity formula. 

\begin{thm} Let $\bf v$ be a bounded solution of \eqref{emain}. Assume furthermore that  $s_*=s_i$ for every $i=1,\cdots,m$. Then 
\begin{equation}
I(R)= \frac{1}{R^{n-2s_*}} \left(\frac{1}{2} \int_{B_R^+}  \sum_{i=1}^{m} y^{a_i} |\nabla v_i|^2 d\mathbf x dy - \int_{B_R\times\{y=0\}} H(\mathbf v) d\mathbf x \right)
\end{equation}
 is a nondecreasing function of $R$ when $H(x_1,\cdots,x_n) \le 0$ for all $(x_1,\cdots,x_n)\in\mathbf R^n$.    
\end{thm}

\noindent \textbf{Proof:} Taking derivative of the function $I(R)$ with respect to $R$ we get 
 \begin{eqnarray*}
I'(R) &=& - \frac{n-2s_*}{2} R^{2s_*-n-1}  \int_{B_R^+} \sum_{i=1}^{m}  y^{a_i} |\nabla v_i|^2 d\mathbf x dy + (n-2s_*)  R^{2s_*-n-1} \int_{B_R \times \{y=0\}} H( \mathbf v) d\mathbf x \\&&+ 
\frac{1}{2}    R^{2s_*-n}  \int_{\partial^+{B_R^+}} \sum_{i=1}^{m}  y^{a_i} |\nabla v_i|^2 d\mathbf x dy - R^{2s_*-n} \int_{\partial {B_R}\times\{y=0\}} H(\mathbf v) d\mathbf x  .
    \end{eqnarray*}
    Multiply (\ref{emain}) with $\mathbf z\cdot\nabla  v_i$ where $\mathbf z=(\mathbf x,y)$ we get 
     \begin{eqnarray*}
0&=& \mathbf z\cdot\nabla v_i \div(y^{a_i} \nabla v_i)
\\&=& \div (y^{a_i}  \nabla v_i \cdot \nabla_{\mathbf z} v_i ) - y^{a_i} \left( |\nabla v_i|^2+\frac{1}{2}\mathbf z\cdot\nabla |\nabla v_i|^2  \right).
      \end{eqnarray*}
Note that $$   y^{a_i} \mathbf z\cdot \nabla |\nabla v_i|^2= \div(y^{a_i}  |\nabla v_i|^2 \mathbf z )- y^{a_i} (n+1)  |\nabla v_i|^2 -a_i y^{a_i} |\nabla v_i|^2  .$$
For each $i=1,\cdots,m$ and $\mathbf z\in\mathbf R^{n+1}_+$,  we have 
\begin{equation}\label{divident}
\div\left( y^{a_i}  \mathbf z\cdot \nabla v_i \nabla v_i   -   \frac{1}{2}  y^{a_i}   |\nabla v_i|^2 \mathbf z \right) + 
\frac{n-2s_i}{2} y^{a_i} |\nabla v_i|^2 =0.
\end{equation}
Integrate over $B_R^+=\{\mathbf X=(\mathbf x,y)\in\mathbf R^{n+1}_+; \ |\mathbf X|<R \}$ where we use the notation $\partial^+B_R^+   =\partial B_R \times \{y>0\}$. Apply the following identities 
 \begin{eqnarray}\label{eq1}
\int_{B_R^+} \div\left( y^{a_i} \nabla v_i \mathbf z\cdot \nabla v_i \right) &=& \int_{B_R\times\{y=0\}} y^{a_i} (-\partial_y v_i) \mathbf z\cdot \nabla v_i + \int_{\partial^+{B_R^+}} y^{a_i} \partial_\nu v_i  \mathbf z\cdot \nabla v_i 
\\ \label{eq2}
\int_{B_R^+} \div\left( y^{a_i} \mathbf z  |\nabla v_i|^2 \right) &=&  \int_{\partial^+{B_R^+}} y^{a_i} \mathbf  z\cdot \mathbf \nu |\nabla v_i|^2
\end{eqnarray}
Note that on $\partial^+ B_R^+$ we have $\mathbf z=R\mathbf \nu$. Using the boundary term in (\ref{emain}) we simplify (\ref{eq1}) and  (\ref{eq2}) as  
 \begin{eqnarray}\label{eq3}
 \int_{B_R^+} \div\left( y^{a_i} \nabla v_i \mathbf z\cdot \nabla v_i \right) &=& \int_{B_R\times\{y=0\}} \partial_{v_i} H(v) \mathbf  x\cdot \nabla_{\mathbf x} v_i +R \int_{\partial^+{B_R^+}} y^{a_i} (\partial_{\mathbf \nu} v_i )^2 d \mathcal{H}^{n} ,
\\ \label{eq4}
\int_{B_R^+} \div\left( y^{a_i} \mathbf z  |\nabla v_i|^2 \right) &=&  R \int_{\partial^+{B_R^+}} y^{a_i}   |\nabla v_i|^2.
\end{eqnarray}
From (\ref{eq3}), (\ref{eq4}) and (\ref{divident}) we get 
 \begin{eqnarray}\label{identity1}
 R \int_{\partial^+{B_R^+}} y^{a_i} (\partial_{\mathbf \nu} v_i )^2 d \mathcal{H}^{n} + \int_{B_R\times\{y=0\}} \partial_{v_i} H(\mathbf v)  \mathbf x\cdot \nabla_{\mathbf x} v_i  \\ \nonumber- \frac{R}{2} \int_{\partial^+{B_R^+}} y^{a_i}   |\nabla v_i|^2 +\frac{n-2s_i}{2} \int_{B_R^+}  y^{a_i} |\nabla v_i|^2 =0. 
 \end{eqnarray}
 Note that 
  \begin{eqnarray}\nonumber
\sum_{i=1}^{m} \int_{B_R} \partial_{v_i} H(\mathbf v)  \mathbf x\cdot \nabla_{\mathbf x} v_i &=& \int_{B_R} \mathbf x\cdot \nabla_{\mathbf x} H(\mathbf v) d\mathbf x  = \int_{B_R}  \div(\mathbf x H(\mathbf v))- n H(\mathbf v) d\mathbf x \\&=& \label{identity2} -n \int_{B_R} H(\mathbf v) d\mathbf x +R \int_{\partial B_R} H(\mathbf v) d \mathcal{H}^{n-1} .
  \end{eqnarray}
Summing (\ref{identity1}) on $i$ and using (\ref{identity2}) we get 
   \begin{eqnarray*}
 R \int_{\partial^+{B_R^+}} \sum_{i=1}^{m} y^{a_i} (\partial_{\mathbf \nu} v_i )^2 d \mathcal{H}^{n} - n \int_{B_R} H(\mathbf v) d\mathbf  x + R \int_{\partial B_R} H(\mathbf v) d \mathcal{H}^{n-1} \\ \nonumber- \frac{R}{2} \int_{\partial^+{B_R^+}} \sum_{i=1}^{m} y^{a_i}   |\nabla v_i|^2 +\sum_{i=1}^{m} \frac{n-2s_i}{2} \int_{B_R^+}  y^{a_i} |\nabla v_i|^2 =0 .
  \end{eqnarray*}
  From this we obtain 
  \begin{eqnarray*}
 \frac{n}{2} \int_{B_R^+}   \sum_{i=1}^{m} y^{a_i} |\nabla v_i|^2  &=&   \sum_{i=1}^{m} s_i \int_{B_R^+}     y^{a_i} |\nabla v_i|^2  + n \int_{B_R} H(\mathbf v) d\mathbf x - R \int_{\partial B_R} H(\mathbf v) d \mathcal{H}^{n-1}  \\&& - R  \int_{\partial^+ B_R^+}  \sum_{i=1}^{m} y^{a_i} (\partial_{\mathbf \nu} v_i )^2 d \mathcal{H}^{n} +  
 \frac{R}{2} \int_{\partial^+{B_R^+}} \sum_{i=1}^{m} y^{a_i}   |\nabla v_i|^2 .
   \end{eqnarray*}
Substitute this in $I'(R)$ we get 
\begin{eqnarray*}
I'(R) R^{n+1-2 s_*} &=& R  \int_{\partial^+ B_R^+}  \sum_{i=1}^{m} y^{a_i} (\partial_\nu v_i )^2 d \mathcal{H}^{n} +   \int_{B_R^+}  \sum_{i=1}^{m} (s_*-s_i)  y^{a_i} |\nabla v_i|^2 - 2 s_* \int_{B_R\times\{y=0\}} H(\mathbf v) d\mathbf x.
   \end{eqnarray*}
Suppose that $s_i=s_*$ for every $i$. Then  $I'(R)\ge 0$ when $H\le 0$. \hfill $ \Box$

\section{Hamiltonian Identities}
In this section we prove a proof for Theorem \ref{hamilton}.  Note that similar Hamiltonian identities for  nonlocal scalar equations  is given by Cabr\'{e}-Sire in \cite{CS} and for the case of local gradient systems, i.e. $\mathbf s=1$ in (\ref{main}), established by Gui in \cite{gui}. 

     Consider a bounded solution $\bf v$ of the extension problem (\ref{emain}) when $s_i=s$ for every $i$. Suppose that $\mathbf v\to\mathbf\alpha$ as $x_n\to\infty$ and $\mathbf v\to\mathbf\beta$ as $x_n\to-\infty$.  

\noindent\textbf{Proof of part (1) of Theorem \ref{hamilton}:} Suppose that $\mathbf v$ is a solution of the extension problem. Define 
\begin{equation}
w(x):=\frac{1}{2} \sum_{i=1}^{m} \int_0^\infty y^{1-2s} \left[  (\partial_x v_i)^2 - (\partial_y v_i)^2 \right] dy.
\end{equation}
Taking the derivative with respect to $x$ we obtain 
\begin{equation}\label{xw}
\partial_x w(x):= \sum_{i=1}^{m} \int_0^\infty y^{1-2s} \left[  \partial_x v_i \partial_{xx} v_i - \partial_y v_i \partial_{xy} v_i \right] dy.
\end{equation}
Note that using the equation (\ref{emain}) we have 
$$ y^a \partial_{xx} v_i + \partial_y \left(   y^a \partial_y v_i  \right)=0. $$
From this and (\ref{xw}) we obtain
\begin{equation}\label{xwi}
\partial_x w(x):= \sum_{i=1}^{m} \int_0^\infty \left[  -\partial_x v_i \partial_y \left(   y^a \partial_y v_i  \right) -  y^a \partial_y v_i \partial_{xy} v_i \right] dy.
\end{equation}
Note that from integration by parts we obtain the following
\begin{equation}
-\int_0^\infty \partial_x v_i \partial_y \left(   y^a \partial_y v_i  \right) dy = \int_0^\infty y^a  \partial_{xy} v_i  \partial_y v_i  dy  -  \lim_{y\to 0} y^a    \partial_{x} v_i \partial_{y} v_i . 
\end{equation}
From this and (\ref{xwi}) we get
 \begin{equation}
\partial_x w(x) = - \sum_{i=1}^{m}  \lim_{y\to 0} y^a    \partial_{x} v_i \partial_{y} v_i .
\end{equation}
From the boundary term in (\ref{main}) we get 
$$\partial_x w(x) = d_s \sum_{i=1}^{m}  \partial_{v_i} H(\mathbf v) \partial_{x} v_i =d_s \partial_{x} \left(H(\mathbf v)\right) .
$$
Therefore, 
$$\partial_x\left(    w(x) - d_s  H(\mathbf v( x,0))  \right)=0.$$
Since $\mathbf v$ is a layer solution, we obtain $$ w(x) -d_s  H(\mathbf v( x,0))+ d_s H(\alpha)=0.$$ If we send $x_n\to-\infty$, then we get that $H(\alpha)=H(\beta)$. 

               \hfill $ \Box$

\noindent\textbf{Proof of part (2) of Theorem \ref{hamilton}:} Suppose that $\mathbf v=\mathbf v(|\mathbf x|,y)$ is radially symmetric in $\mathbf x$. For every $i$ and $r=|\mathbf x|$ we have
 \begin{eqnarray}\label{emainrad}
 \left\{ \begin{array}{lcl}
\hfill \partial_{rr} v_i +\frac{n-1}{r} \partial_r v_i + \partial_{yy} v_i +\frac{a_i}{y} \partial_y v_i&=& 0   \ \ \text{in}\ \  (0,\infty)\times(0,\infty),\\   
\hfill -\lim_{y\to0}y^{1-2s_i} \partial_{y} v_i&=&  d_{s_i} \partial_{v_i}H(\mathbf v)   \ \ \text{in}\ \ (0,\infty)\times(y=0).
\end{array}\right.
  \end{eqnarray}
Set $a=1-2s$.  Define the following function of $r$, 
  \begin{equation}
  w(r):= \sum_{i=1}^{m} \int_0^\infty   \frac{y^{1-2s}}{2} \left[  (\partial_r v_i)^2 - (\partial_y v_i)^2 \right] dy .
  \end{equation}
 Taking derivative of the above with respect to $r$ we get 
 \begin{equation}
  \partial_r w(r):= \sum_{i=1}^{m} \int_0^\infty   y^{1-2s} \left[  \partial_r v_i  \partial_{rr} v_i - \partial_y v_i \partial_{ry} v_i \right] dy .
  \end{equation}
 Note that using the equation (\ref{emainrad}) we have 
 $$ \partial_{rr} v_i = -\frac{n-1}{r} \partial_r v_i - \partial_{yy} v_i - \frac{a}{y} \partial_y v_i$$
 Substituting this in $\partial_r w(r)$ we obtain
 \begin{eqnarray}\label{wr}
\nonumber \partial_r w(r) &=&   -\frac{n-1}{r}  \sum_{i=1}^{m} \int_0^\infty   y^{a} (\partial_r v_i)^2 dy -   \sum_{i=1}^{m} \int_0^\infty   y^{a}  \partial_{yy} v_i \partial_{r} v_i dy \\&& - a \sum_{i=1}^{m} \int_0^\infty   y^{a-1} \partial_{r} v_i \partial_{y} v_i dy -  \sum_{i=1}^{m} \int_0^\infty   y^{a}  \partial_{y} v_i  \partial_{ry} v_i .
  \end{eqnarray}
  Applying integration by parts we obtain 
  \begin{eqnarray*}
  \int_0^\infty   y^{a}  \partial_{y} v_i  \partial_{ry} v_i = -  \int_0^\infty  \partial_{y}\left(y^{a}  \partial_{y} v_i\right)  \partial_{r} v_i dy + \lim_{y\to 0} y^a  \partial_{r} v_i \partial_{y} v_i.
  \end{eqnarray*}
that gives us 
  \begin{eqnarray*}
- \lim_{y\to 0} y^a  \partial_{r} v_i \partial_{y} v_i  &=& -  \sum_{i=1}^{m} \int_0^\infty   y^{a}  \partial_{yy} v_i \partial_{r} v_i dy \\&& - a \sum_{i=1}^{m} \int_0^\infty   y^{a-1} \partial_{r} v_i \partial_{y} v_i dy -  \sum_{i=1}^{m} \int_0^\infty   y^{a}  \partial_{y} v_i  \partial_{ry} v_i .
   \end{eqnarray*}
From this and (\ref{wr}) we get 
 \begin{equation}\label{partialrw} 
 \partial_{r} w(r) =  -\frac{n-1}{r}  \sum_{i=1}^{m} \int_0^\infty   y^{a} (\partial_r v_i)^2 dy - \lim_{y\to 0} y^a  \sum_{i=1}^{m}  \partial_{r} v_i \partial_{y} v_i .
 \end{equation}
We now apply the boundary term in (\ref{emainrad}) to get 
$$- \lim_{y\to 0} y^a  \sum_{i=1}^{m}  \partial_{r} v_i \partial_{y} v_i = \frac{1}{d_s} \sum_{i=1}^{m} \partial_{v_i} H(\mathbf v) \partial_{r} v_i  = d_s  \partial_{r} \left(H(\mathbf v)\right). $$
 From this and (\ref{partialrw}) we get the following 
$$ \partial_{r} \left( w(r) - \frac{1}{d_s}  \partial_{r} \left(H(\mathbf v)\right) \right)=   -\frac{n-1}{r}  \sum_{i=1}^{m} \int_0^\infty   y^{a} (\partial_r v_i)^2 dy <0.$$ 
  This completes the proof. 

               \hfill $ \Box$

\section{Structure of the nonlinearity $H(\mathbf u)$ for radial solutions}
This section is devoted to the proof of Theorem \ref{radialProp}.  The proof follows the one in the scalar case \cite{CS1}.  Consider  the translated (or slided) functions $\mathbf v^t(\mathbf x,y)= \mathbf v(x_1+t, \cdots, x_n+t,y)$. It is straightforward to see that these functions solve 
the extension problem (\ref{emain}).  By the H\"{o}lder estimates in Lemma \ref{regv}, the translated solutions converge locally uniformly and up to subsequences, to a solution of the same problem (\ref{emain}). From the assumption  $\lim_{|\mathbf x|\to \infty} v(|\mathbf x|,0)=0$ such limit is identically constant. From this for each $i$ and $|\mathbf x|\to 0$  we conclude that  
$$||v_i||_{L^{\infty}(B_R^+(\mathbf x,0))} + || \nabla_{\mathbf x} v_i||_{L^{\infty}(B_R^+(\mathbf x,0))} + ||y^{a_i} \partial_y v_i||_{L^{\infty}(B_R^+(\mathbf x,0))}\to 0 $$
that implies $$\partial_{v_i} H(0)=0.$$
We now apply the monotonicity formula (\ref{radialmono}) when $r\to\infty$ and $r=0$ to get  \begin{equation}\label{radialmono}
   -\sum_{i=1}^{m} \int_0^\infty y^{1-2s}  (\partial_y v_i)^2 dy - 2 d_s H(\mathbf v(0,0)) \ge -2 d_s  H(0)
  \end{equation}
Therefore, $H(\mathbf v(0,0)) \le  H(0)$. 

Now differentiate the equation (\ref{emain}) with respect to $r=|\mathbf x|$ to get 
$$ \partial_{rr} v_i +\frac{n-1}{r} \partial_r v_i + \partial_{yy} v_i +\frac{a_i}{y} \partial_y v_i= 0 \ \ (0,\infty)\times(0,\infty).$$
Set $\phi_i:=-u_r>0$ in $(\mathbf R^n\setminus\{0\}) \times (0,\infty)$.  Straightforward calculations show that 
 \begin{eqnarray}
 \left\{ \begin{array}{lcl}
\hfill \div(y^{a_i} \nabla \phi_i) &=& \frac{n-1}{|\mathbf x|^2} y^{a_i} \phi_i    \ \ \text{in}\ \  (\mathbf R^n\setminus\{0\}) \times (0,\infty),\\   
\hfill -\lim_{y\to0}y^{1-2s_i} \partial_{y} \phi_i&=&  d_{s_i} \sum_{j=1}^m \partial_{v_i v_j}H(\mathbf v) \phi_j  \ \ \text{in}\ \ (\mathbf R^n\setminus\{0\}) \times(y=0)
\end{array}\right.
  \end{eqnarray}
  For any $\mathbf x_0\in\mathbf R^n$, define the shift functions $\phi_i^{\mathbf x_0}:=\phi_i(\mathbf x-\mathbf x_0,y)>0$ in $(\mathbf R^n\setminus\{\mathbf x_0\}) \times (0,\infty)$ that satisfies
   \begin{eqnarray}\label{phix0}
   \left\{ \begin{array}{lcl}
\hfill \div(y^{a_i} \nabla \phi_i^{\mathbf x_0}) &=& \frac{n-1}{|\mathbf x - \mathbf x_0|^2} y^{a_i} \phi_i^{\mathbf x_0}   \ \ \text{in}\ \  (\mathbf R^n\setminus\{0\}) \times (0,\infty),\\   
\hfill -\lim_{y\to0}y^{1-2s_i} \partial_{y} \phi_i^{\mathbf x_0} &=&  d_{s_i} \sum_{j=1}^m \partial_{v_i v_j}H(\mathbf v^{\mathbf x_0}) \phi_j^{\mathbf x_0}  \ \ \text{in}\ \ (\mathbf R^n\setminus\{0\}) \times(y=0)
\end{array}\right.
  \end{eqnarray}
  where $v_i^{\mathbf x_0}:=v_i(\mathbf x-\mathbf x_0,y)$ for every $i$.   Consider the test function $\zeta\in C^1(\mathbf R^n)$ such that it vanishes on $\partial B_R \times [0,R)$ and $B_R\times \{y=R\}$. Then for any $|\mathbf x_0|>R$ multiply the boundary equation of  (\ref{phix0}) with $\frac{\zeta^2}{\phi_i^{\mathbf x_0}}$ and integrate over $B_R$ to obtain
     \begin{eqnarray}\label{}
 d_{s_i} \sum_{j=1}^m \int_{B_R} \partial_{v_i v_j}H(\mathbf v^{\mathbf x_0}) \phi_j^{\mathbf x_0} \frac{\zeta^2}{\phi_i^{\mathbf x_0}} &=& \int_{C_R} \div(y^{a_i} \nabla \phi_i^{\mathbf x_0}) \frac{\zeta^2}{\phi^{\mathbf x_0}} +y^{a_i} \nabla \phi_i^{\mathbf x_0} \cdot \nabla \left(\frac{\zeta^2}{\phi_i^{\mathbf x_0}}\right) 
 \\&=&  \int_{C_R} \frac{n-1}{|\mathbf x - \mathbf x_0|^2} y^{a_i} \zeta^2 + y^{a_i} \nabla \phi_i^{\mathbf x_0} \cdot \nabla \left(\frac{\zeta^2}{\phi_i^{\mathbf x_0}}\right) 
  \\&=&  \int_{C_R}   \frac{n-1}{|\mathbf x - \mathbf x_0|^2} y^{a_i} \zeta^2 + y^{a_i} \left[ 2\frac{\zeta}{ \mathbf x_0 }  \nabla \zeta \cdot \nabla \phi_i^{\mathbf x_0} - \zeta^2 \frac{  |\nabla \phi_i^{\mathbf x_0} |^2  }{ |\phi_i^{\mathbf x_0} |^2  }\right]
\\&  \le & \int_{C_R}    \frac{n-1}{|\mathbf x - \mathbf x_0|^2} y^{a_i} \zeta^2 + y^{a_i} |\nabla \zeta|^2
     \end{eqnarray}
     Sending $|\mathbf x_0|\to \infty$ then 
\begin{equation}\label{H0}
   \sum_{j=1}^m \partial_{v_i v_j}H(0) \int_{B_R}  \phi_j^{\mathbf x_0} \frac{\zeta^2}{\phi_i^{\mathbf x_0}} \le \frac{1}{d_{s_i}}  \int_{C_R}  y^{a_i} |\nabla \zeta|^2.
\end{equation}
Taking the sum on $i$ we get 
     \begin{eqnarray*}\label{}
   \sum_{i,j=1}^m \partial_{v_i v_j}H(0)  \frac{\phi_j^{\mathbf x_0}}{\phi_i^{\mathbf x_0}}   \zeta^2 &=& \sum_{i=1}^m \partial_{v_i v_i}H(0)  \zeta^2      + \sum_{i<j} \partial_{v_i v_j}H(0)   \left(\frac{ \phi_j^{\mathbf x_0} }{\phi_i^{\mathbf x_0}}  +   \frac{ \phi_i^{\mathbf x_0} }{\phi_j^{\mathbf x_0}}\right  )  \zeta^2
   \\&\ge&  \sum_{i,j=1}^m  \partial_{v_i v_j}H(0) \zeta^2.
     \end{eqnarray*}
From this and (\ref{H0}), we obtain 
\begin{equation}\label{H01}
 \sum_{i,j=1}^m  \partial_{v_i v_j}H(0) \le   \frac{ \sum_{i=1}^m   \frac{1}{d_{s_i}}  \int_{C_R}  y^{a_i} |\nabla \zeta|^2 d\mathbf x dy   }{\int_{B_R} \zeta^2(\mathbf x,0) d\mathbf x }.
\end{equation}
Suppose that $\eta_R>0$ is the first eigenfunction of the operator $-\Delta_{\mathbf x} $ in $B_R$ with Dirichlet boundary conditions.  Now define the smooth function $\chi_R(y): \mathbf R^+\to  [0,1]$ that satisfies  $\chi_R(y)=1$ for  $y\in [0,R/2]$ and $\chi_R(y)=0$ for $y\in [R,\infty)$ and $||\chi_R'||_{L^\infty(R/2,R)}<C/R$. Now set $\zeta_R(\mathbf x,y):= \eta_R(\mathbf x) \chi_R(y)$ and test (\ref{H01}) on $\zeta_R$. Note that $\zeta_R$ is zero on $\partial B_R \times [0,R)$ and $B_R\times \{y=R\}$.  For the left-hand side we have
     \begin{eqnarray}\label{nzeta}
 \int_{C_R}  y^{a_i} |\nabla \zeta|^2 &=& \int_{C_R}  y^{a_i}  \left[  |\nabla _{\mathbf x}\eta_R|^2 \chi_R^2 + \eta_R^2 |\chi'_R|^2  \right]
 \\&=& \nonumber \left[ \lambda_R  \int_0^R y^{a_i} \chi_R^2 dy + \int_0^R y^{a_i} |\chi'_R|^2 dy  \right] \int_{B_R} \eta_R^2.
     \end{eqnarray}
Note that $\lambda_R=\frac{C(n)}{R^2}$ is the first eigenvalue of the operator $-\Delta_{\mathbf x} $ in $B_R$ with Dirichlet boundary conditions and $\eta_R(\mathbf x)=\zeta_R(\mathbf x,0)$ because $\chi_R(0)=1$. From this and (\ref{nzeta}) we have 
\begin{eqnarray*}\label{}
 \int_{C_R}  y^{a_i} |\nabla \zeta|^2 d\mathbf x dy&=&  \nonumber \left[ C(n) R^{-2}  \int_0^R y^{a_i} \chi_R^2 dy + \int_0^R y^{a_i} |\chi'_R|^2 dy  \right] \int_{B_R} \zeta_R^2(\mathbf x,0)d\mathbf x
 \\&\le& C R^{-2} \int_0^R y^{a_i}  dy \int_{B_R} \zeta_R^2 (\mathbf x,0)d\mathbf x
  \\&\le&  C (1-s_i)^{-1} R^{-2s_i} \int_{B_R} \zeta_R^2(\mathbf x,0)d\mathbf x.
     \end{eqnarray*}
Therefore, for every $i$ we have 
$$  \frac{1}{d_{s_i}}  \int_{C_R}  y^{a_i} |\nabla \zeta|^2 d\mathbf x dy \le   C R^{-2s_i} \int_{B_R} \zeta_R^2(\mathbf x,0)d\mathbf x $$
where $C=C(n,s_i)$ is independent of $R$. Sending $R\to\infty$ and using (\ref{H01}) we get $$  \sum_{i,j=1}^m  \partial_{v_i v_j}H(0) \le 0.$$

\section{Liouville theorems}
In this section we assume that the right-hand side of the nonlinearity in system (\ref{main}) has a sign, i.e. $\nabla H(u)\ge 0$. For such systems, we can provide a Liouville theorem for general nonlinearities as Theorem \ref{lioupositive}.  The proof of the theorem follows from several lemmata.  For the scalar case, $m=1$, similar methods are developed by Dupaigne-Sire in \cite{DS}. Also, similar Liouville theorems are given for the scalar local problem, $s=m=1$, by Dupaigne-Farina in \cite{DF}. 

\begin{lemma}\label{vibound}
Suppose that $\mathbf v=(v_i)_i$ is a bounded solution of (\ref{emain}) where $\nabla H \ge0$. Then, there exists a constant $C>0$ such that for all $R>1$ and any $i=1,\cdots,m$
\begin{equation}\label{EnHpos}
\int_{B_R^+} y^{1-2s_i} |\nabla v_i|^2 \le C R^{n-2s_i}
\end{equation}
\end{lemma}

 \noindent\textbf{Proof:} Since $v_i$ is bounded, multiply the equation with $(v_i-||v_i||_{\infty})  \phi_R$ where $\phi_R$ is a positive test function.  This gives us  
\begin{equation}\label{EnHpos}
\int_{\mathbf R^{n+1}_+}  (v_i-||v_i||_{\infty})  \phi_R \div(y^{a_i} \nabla v_i)=0.
\end{equation}
Therefore, 
\begin{eqnarray}\label{EnHpos1}
\nonumber\int_{\mathbf R^{n+1}_+} y^{a_i}  \nabla[ (v_i-||v_i||_{\infty})  \phi_R]\cdot\nabla v_i  &=&\int_{\partial\mathbf R^{n+1}_+}  y^{a_i}  (v_i-||v_i||_{\infty})    (-\partial_y v_i) \phi_R
\\&=&\int_{\partial\mathbf R^{n+1}_+}  y^{a_i}  (v_i-||v_i||_{\infty}) \partial_{v_i} H(\mathbf v(\mathbf x,0)) \phi_R \le 0,
\end{eqnarray}
where the assumption $\nabla H\ge 0$ is used.  From this we get 
\begin{eqnarray*}
\int_{\mathbf R^{n+1}_+} y^{a_i} |\nabla v_i|^2 \phi_R &\le& - \int_{\mathbf R^{n+1}_+} y^{a_i} (v_i-||v_i||_{\infty})   \nabla \phi_R\cdot \nabla v_i \\&=& \frac{1}{2}  \int_{\mathbf R^{n+1}_+} y^{a_i} \nabla (v_i-||v_i||_{\infty})^2\cdot \nabla \phi_R \\&=&- \frac{1}{2}  \int_{\mathbf R^{n+1}_+} (v_i-||v_i||_{\infty})^2 \div( y^{a_i}\nabla \phi_R) + \frac{1}{2} \int_{\partial\mathbf R^{n+1}_+} (v_i-||v_i||_{\infty})^2 (-  y^{a_i} \partial_y \phi_R)
\\&=:& I+J
\end{eqnarray*}
Using the fact that $v_i$ is bounded and applying an appropriate test function $\phi_R$ we have 
\begin{eqnarray*}
I&=&- \frac{1}{2}  \int_{\mathbf R^{n+1}_+} (v_i-||v_i||_{\infty})^2 \left( y^{a_i} \Delta_{\mathbf x} \phi_R +a_i y^{a_i-1} \partial_y \phi_R   \right)
\\&\le& C  \int_{B_{2R}^+} \left(\frac{y^{a_i}}{R^2} +\frac{y^{a_i-1}}{R} \right)  d\mathbf x dy
\\&\le& C R^{n} \int_{R}^{2R}  \left(\frac{y^{a_i}}{R^2} +\frac{y^{a_i-1}}{R} \right)  dy
\\&\le& C R^{n-2s_i}
\end{eqnarray*}
where the constant $C$ is independent of $R$. Note that we have used the fact that the test function $\phi_R$ can be chosen such that $|\partial_y \phi_R |\le CR^{-1}$ and $|\Delta_{\mathbf x} \phi_R| \le C R^{-2}$ on $B_R^+$. Similarly, using these properties we have  $ |y^{a_i}  \partial_y \phi_R |\le C R^{-2s_i}$.  This yields 
$$ J \le C R^{n-2s_i}.$$
This finishes the proof. 

               \hfill $ \Box$
               
The next lemma is a characterization of bounded pointwise-stable solutions  of (\ref{emain}) when $n=1$. 

\begin{lemma}\label{ux}
Suppose that $\mathbf v=(v_i)_i$ is a bounded pointwise-stable solution of (\ref{emain}) when $n=1$. Then, for every $1\le i \le m$   either each $\partial_x v_i$ vanishes  in $\overline{\mathbf{R}^{2}_+}$ or it does not change sign in $\overline{\mathbf{R}^{2}_+}$. 
\end{lemma}
 \noindent\textbf{Proof:} The function $v_i$ is a pointwise-stable solution of  
\begin{eqnarray}
 \left\{ \begin{array}{lcl}
\hfill \div(y^{a_i} \nabla v_i)&=& 0   \ \ \text{in}\ \ \mathbf{R}_+^{2},\\   
\hfill -\lim_{y\to0}y^{a_i} \partial_{y} v_i&=& d_{s_i} \partial_{v_i}H(\mathbf v)   \ \ \text{in}\ \ \partial\mathbf{R}_+^{2}.
\end{array}\right.
  \end{eqnarray}
  Therefore, there exits $\mathbf\phi=(\phi_i)_i$ such that for each $i$ the function $\phi_i$ does not change sign and it satisfies  
  \begin{eqnarray}\label{}
 \left\{ \begin{array}{lcl}
\hfill \div(y^{a_i} \nabla \phi_i)&=& 0   \ \ \text{in}\ \ \mathbf{R}_+^{2},\\   
\hfill -\lim_{y\to0}y^{a_i} \partial_{y} \phi_i&=& \sum_{i=1}^m d_{s_i} \partial_{v_i v_j}H(\mathbf v) \phi_j   \ \ \text{in}\ \ \partial\mathbf{R}_+^{2}. 
\end{array}\right.
  \end{eqnarray}
   For each $i$, set $\sigma_i=\frac{\partial_x v_i}{\phi_i}$. Then 
   \begin{eqnarray*}
\div(y^{a_i} \phi_i^2 \nabla \sigma_i) &=& \div(y^{a_i} \phi_i^2 \nabla v_i ) -  \div(y^{a_i} \partial_{x} v_i \nabla\phi_i) \\&=&  \div(y^{a_i} \nabla v_i )\phi_i + y^{a_i} \nabla \phi_i \cdot \nabla \partial_x v_i -  \div(y^{a_i} \nabla \phi_i )\partial_x v_i - y^{a_i} \nabla \phi_i \cdot \nabla \partial_x v_i \\&=&0. 
      \end{eqnarray*}
It is straightforward to see that $-\lim_{y\to0}y^{a_i} \partial_{y} \sigma_i=0$ on $\partial\mathbf 
R^{2}_+$. Now multiply the equation (\ref{emain}) by $v_i\phi_R$ where $\phi_R$ is a test function and integrate over $\mathbf R^{2}_+$. Then, 
$$  \sum_{i=1}^{m} \int_{\mathbf R^{2}_+} y^{a_i} |\nabla v_i|^2 \phi_R^2+ 2  \sum_{i=1}^{m}\int_{\mathbf R^{2}_+} y^{a_i} \nabla \phi_R\cdot\nabla v_i \phi_R v_i= \sum_{i=1}^{m} \int_{\partial\mathbf R^{2}_+} d_{s_i} \partial_{v_i} H(\mathbf v( x,0)) v_i \phi_R^2.$$
Using Cauchy-Schwarz inequality we have
$$ 2 \int_{\mathbf R^{2}_+} y^{a_i} \nabla \phi_R\cdot\nabla v_i \phi_R v_i \le \frac{1}{2} \int_{\mathbf R^{2}_+} y^{a_i} |\nabla v_i|^2 \phi_R^2 + 2 \int_{\mathbf R^{2}_+} y^{a_i} |\nabla \phi_R|^2 v_i^2.$$
From this we obtain
$$ \sum_{i=1}^{m} \int_{\mathbf R^{2}_+} y^{a_i} |\nabla v_i|^2 \phi_R^2 \le C \sum_{i=1}^{m} \int_{\mathbf R^{2}_+} y^{a_i} |\nabla \phi_R|^2 v_i^2+ C \sum_{i=1}^{m} \int_{\mathbf R} d_{s_i} \partial_{v_i} H(\mathbf v(\mathbf x,0)) v_i \phi_R^2.
$$
From the boundedness of $\mathbf v$ and the fact that $\partial_{v_i} H$ is regular enough, we obtain the following bound by choosing the standard test function, 
$$ \sum_{i=1}^{m} \int_{B_R^+} y^{a_i} |\nabla v_i|^2 \phi_R^2 \le C (R^{1-2s_i} + R) \le C R .$$
 Note that $(\sigma_i\phi_i)^2=(\partial_x v_i)^2$ and therefore 
$$ \sum_{i=1}^{m}  \int_{B_R^+}  y^{a_i} (\sigma_i \phi_i)^2 \le C   \sum_{i=1}^{m}  \int_{B_R^+}  y^{a_i} |\nabla v_i|^2 \le C R. $$
Applying Liouville Theorem \ref{liouville}  we obtain that each $\sigma_i$ is contact that is there exists constant $ C_i$ such that $\partial_x v_i=C_i \phi_i$ where $\phi_i$ does not change sign. This finishes the proof.   

 \hfill $ \Box$

 \noindent\textbf{Proof of Theorem \ref{lioupositive}:} The proof follows from Lemma \ref{ux} and Lemma \ref{vibound}.  Suppose that $n \le 2(1+s_*)$ and therefore for every $1\le i\le m$ we have $n-2 s_i\le n-2s_*$. From Lemma \ref{vibound} and Theorem \ref{thsymv} we conclude that $\mathbf v$ is a bounded pointwise-stable solution of (\ref{emain}) for $n=1$. Then, Lemma \ref{ux} implies that for every $1\le i \le m$   either each $\partial_x v_i$ vanishes  in $\overline{\mathbf{R}^{2}_+}$ or it does not change sign in $\overline{\mathbf{R}^{2}_+}$. If $\partial_x v_i$ vanishes in $\overline{\mathbf{R}^{2}_+}$ then there are constants $C_i,\hat C_i$ such that $v_i(x,y)=\hat C_i y^{1-s_i}+C_i$ and by the boundedness of $v_i$ we have $\hat C_i=0$.  So, $v_i\equiv C_i$. 
 
 Now suppose that $v_i(x,y)$ is strictly monotone in $x$. From boundedness of $v_i$, we conclude that $v_i(x,0)$ has limits for $x\to\pm\infty$. Let's define $\lim_{x\to\infty} v_i=\alpha_i$ and $\lim_{x\to-\infty} v_i=\beta_i$. Since  $v_i(x,0)$ is strictly monotone in $x$, for every $1\le i\le m$ we have $\beta_i<\alpha_i$.  Therefore, $\mathbf v(x,y)=(v_i(x,y))_i$ is a bounded layer solution of (\ref{emain}) for $n=1$ where $\lim_{x\to\infty} \mathbf v=\alpha=(\alpha_i)_i$ and $\lim_{x\to-\infty} \mathbf v=\beta=(\beta_i)_i$. Applying Theorem \ref{hamilton} we conclude that for any $x\in \mathbf R$ the following Hamiltonian identity holds.
\begin{equation}\label{hamiltonian1}
 \sum_{i=1}^{m} \int_0^\infty y^{1-2s} \left[  (\partial_x v_i)^2 - (\partial_y v_i)^2 \right] dy=2d_s \left[  H(\mathbf v(x,0)) - H(\mathbf \alpha) \right]. 
\end{equation}
Now sending $x\to-\infty$ and using the boundedness of $\mathbf v$ and also applying Lemma \ref{decayv}, from the Hamiltonian identity (\ref{hamiltonian1}) we obtain $H(\alpha)=H(\beta)$.   From the assumptions we already know that $\nabla H \ge 0$ and $H$ is not identically constant.  This is in fact in contradiction with the mean value theorem that states there is a $t\in (0,1)$ such that $$0=H(\alpha)-H(\beta)=(\alpha-\beta)\cdot\nabla H(t(\alpha-\beta)+\beta)>0.$$ 
This finishes the proof. 
  
 \hfill $ \Box$

{\bf Acknowledgment} Both authors appreciate Professor Juncheng Wei for his invitation to the University of British Columbia, where this work was initiated, and for his supports and fruitful discussions.

\bibliographystyle{plain}
\bibliography{biblio}

       \end{document}